\input ./style/arxiv-general.cfg
\documentclass[aos,MSNbibl,seceqn,nameyear,dvips]{arximspdf}
\makeatletter
   \@ifpackageloaded{graphicx}{}{\usepackage{graphicx}}
\makeatother
\usepackage{dcolumn}
\usepackage{multirow}

\doi{10.1214/15-AOS1380}% Updated by VTEXPTS2LaTeX.exe, 20.10.2015 14:42
\volume{44}
\issue{2}
\pubyear{2016}
\firstpage{629}
\lastpage{659}
\docsubty{FLA}

\makeatletter
\newcolumntype{d}[1]{D{.}{.}{#1}}
\newcommand{\lleft}{\left}
\newcommand{\rrvert}{\vert}
\newcommand{\rright}{\right}
\newcommand{\rrVert}{\Vert}
\newcommand{\llvert}{\vert}
\newcommand{\llVert}{\Vert}
\newcommand{\citepp}[1]{[\citet{#1}]}
\newcommand{\eqref}[1]{(\ref{#1})}
\newtheorem{theorem}{Theorem}[section]
\newtheorem{lemma}[theorem]{Lemma}
\newtheorem{proposition}[theorem]{Proposition}
\newtheorem{corollary}[theorem]{Corollary}
\makeatother

\begin{document}
\begin{frontmatter}

\title{Global solutions to folded concave penalized nonconvex learning}
\runtitle{Global solutions to nonconvex learning}

\begin{aug}
% Corresponding author: Hongcheng Liu - hql5143@psu.edu% Updated by
%VTEXPTS2LaTeX.exe, 20.10.2015 14:42
\author[A]{\fnms{Hongcheng}~\snm{Liu}\corref{}\thanksref{T2}\ead[label=e1]{hql5143liu@gmail.com}},
\author[A]{\fnms{Tao}~\snm{Yao}\thanksref{T2}\ead[label=e2]{taoyao@psu.edu}}
\and
\author[B]{\fnms{Runze}~\snm{Li}\thanksref{T4}\ead[label=e3]{rzli@psu.edu}}
\runauthor{H. Liu, T. Yao and R. Li}
\affiliation{Pennsylvania State University}
%\dedicated{}
\address[A]{H. Liu\\
T. Yao\\
Department of Industrial\\
\quad and Manufacturing Engineering\\
Pennsylvania State University\\
University Park, Pennsylvania 16802\\
USA\\
\printead{e1}\\
\phantom{E-mail: }\printead*{e2}}
\address[B]{R. Li\\
Department of Statistics\\
Pennsylvania State University\\
University Park, Pennsylvania 16802\\
USA\\
\printead{e3}}
\end{aug}
\thankstext{T2}{Supported by Penn State Grace Woodward Collaborative
Engineering/Medicine Research grant,
NSF Grant CMMI 1300638, Marcus PSU-Technion
Partnership grant and Mid-Atlantic University Transportation Centers grant.}
\thankstext{T4}{Support by NSF Grant DMS-15-12422 and
National Institute of Health Grants P50 DA036107 and P50 DA039838.}

% HISTORY:
%
\received{\smonth{7} \syear{2014}}% Updated by VTEXPTS2LaTeX.exe,
%20.10.2015 14:42
%
\revised{\smonth{6} \syear{2015}}% Updated by VTEXPTS2LaTeX.exe,
%20.10.2015 14:42

% ABSTRACT
%
\begin{abstract}
This paper is concerned with solving nonconvex learning problems with
folded concave penalty.
Despite that their global solutions entail desirable statistical
properties, they lack optimization
techniques that guarantee global optimality in a general setting. In this
paper, we show that a class of nonconvex learning problems are
equivalent to
general quadratic programs. This equivalence facilitates us in developing
mixed integer linear programming reformulations, which admit finite algorithms
that find a provably global optimal solution. We refer to this
reformulation-based technique as the mixed integer programming-based global
optimization (MIPGO). To our knowledge, this is the first global optimization
scheme with a theoretical guarantee for folded concave penalized nonconvex
learning with the SCAD penalty [\textit{J.~Amer. Statist. Assoc.} \textbf{96} (2001)
1348--1360] and the MCP penalty [\textit{Ann. Statist.} \textbf{38} (2001) 894--942].
Numerical results indicate a significant outperformance of
MIPGO over
the state-of-the-art solution scheme, local linear approximation and other
alternative solution techniques in literature in terms of solution quality.
\end{abstract}

% KEYWORDS
\begin{keyword}[class=AMS]
\kwd[Primary ]{62J05}
%\kwd{}
\kwd[; secondary ]{62J07}
\end{keyword}
\begin{keyword}
\kwd{Folded concave penalties}
\kwd{global optimization}
\kwd{high-dimensional statistical learning}
\kwd{MCP}
\kwd{nonconvex quadratic programming}
\kwd{SCAD}
\kwd{sparse recovery}
\end{keyword}
\end{frontmatter}

%s1 #&#
\section{Introduction}

Sparse recovery is of great interest in high-dimensional statistical
learning. Among the most
investigated sparse recovery techniques are LASSO and the nonconvex
penalty methods, especially
folded concave penalty techniques [see \citet{FanandLv2011}, for a general
definition]. Although LASSO is a
popular tool primarily because its global optimal solution is
efficiently computable, recent theoretical and
numerical studies reveal that this technique requires a critical
irrepresentable condition to
ensure statistical performance. In comparison, the folded concave
penalty methods require less theoretical regularity
and entail better statistical properties \citepp
{Zou2006,Meinshausen,Fanetal2012}.
%
%
%\citepp{ZhaoandYu2006,Zou2006,Meinshausen,FanandLi2001,FanandPeng2004,Zhang2010,FanandLv2011,Fanetal2012}.
In particular, \citet{ZhangandZhang2012} showed that the global
solutions to the folded concave penalized learning problems lead to a
desirable recovery performance. However, these penalties cause the
learning problems to be nonconvex and render the local solutions to be
nonunique in general.

Current solution schemes in literature focus on solving a nonconvex
learning problem locally.
\citet{FanandLi2001} proposed a local quadratic approximation (LQA)
method, which was further analyzed by using majorization minimization
algorithm-based techniques in \citet{HunterandLi2005}. \citet
{Mazumder2011} and \citet{BrehenyandHuang2011} developed different
versions of coordinate descent algorithms. \citet{ZouandLi2008} proposed
a local linear approximation (LLA) algorithm
and \citet{Zhang2010} proposed a PLUS algorithm. \citet{KimCCCP}
developed the ConCave Convex procedure (CCCP).
To justify the use of local algorithms, conditions were imposed for the
uniqueness of a local solution \citepp{Zhang2010,ZhangandZhang2012}; or,
even if multiple local minima exist, the
strong oracle property can be attained by LLA with wisely (but fairly
efficiently) chosen initial solutions \citepp{Fanetal2012}. \citet
{HuangandZhang2012} showed that a multistage framework that subsumes
the LLA can improve the solution quality stage by stage under some
conditions. \citet{Wangetal2013a} proved that calibrated CCCP
produces a
consistent solution path which contains the oracle estimator with
probability approaching one. \citet{LohandWainwright} established
conditions for all local optima to lie within statistical
precision of the true parameter vector, and proposed to employ the gradient
method for composite objective function minimization by \citet
{Nesterov2007} to
solve for one of the local solutions. \citet{Wangetal2013}
incorporated the
gradient method by \citet{Nesterov2007} into a novel approximate
regularization path following algorithm, which was shown to converge linearly
to a solution with an oracle statistical property. Nonetheless, none of
the above algorithms theoretically ensure global optimality. %there has
%been scarce discussions on solution schemes that theoretically
%guarantees global optimality.

%Nonetheless, the convergence property and the statistical performance
%analysis in the aforementioned works can be contingent on certain
%imposed conditions, some of which are artificially imposed. More
%specifically, \citet{Fanetal2012} involve a restricted eigenvalue
%(RE)
%condition necessary for its convergence (to the oracle solution) with
%lower bounded probability. Both \citet{LohandWainwright} and
%\citet{Wangetal2013} require restricted strong convexity (RSC) to
%ensure convergence. Without these conditions, the local solution
%methods may not have a theoretical guarantee in solution quality. This
%will be demonstrated in our numerical comparison in Section~
%\ref{section nesterov} and \ref{section LLA}.

In this paper, we seek to solve folded concave penalized nonconvex
learning problems in a direct and generic way: to derive a reasonably
efficient solution scheme with a provable guarantee on global
optimality. Denote by $n$ the
sample size, and by $d$ the problem dimension. Then the folded
concave penalized learning problem of our discussion is formulated as following:
%e1.1 #&#
%
\begin{equation}
\min_{\beta\in\Lambda}\mathcal L(\beta):= \mathbb L(\beta)+n\sum
_{i=1}^dP_\lambda\bigl(\llvert
\beta_i\rrvert \bigr), \label{LP folded concave}
\end{equation}
where $P_\lambda( \cdot): \mathbb R\rightarrow\mathbb R$ is a
penalty function with tuning parameter $\lambda$. Our
proposed procedure is directly applicable for settings allowing $\beta
_i$ to have different $\lambda_i$ or different
penalty. For ease of presentation and without loss of generality, we
assume $P_\lambda( \cdot)$ is the same
for all coefficients. Function $\mathbb L( \cdot):\mathbb
R^d\rightarrow\mathbb R$ is defined as a quadratic function,
$\mathbb L(\beta):= \frac{1}{2}\beta^\top Q\beta+q^\top\beta$,
which is an abstract
representation of a proper (quadratic) statistical loss function with $
Q\in\mathbb R^{d\times
d}$ and $q\in\mathbb R^d$ denoting matrices from data samples.
Denote by $\Lambda:= \{\beta\in\mathbb R^d: \mathcal A^\top\beta
\leq\mathbf b\}$ the
feasible region defined by a set of linear constraints with $\mathcal
A\in\mathbb
R^{d\times m}$ and $\mathbf b\in\mathbb R^m$ for some proper $m:
0\leq m<d$. Assume throughout the paper that $Q$ is symmetric, $
\mathcal A$ is full rank and $\Lambda$ is nonempty.
Notice that under this
assumption, the loss function does not have to be convex. We instead stipulate
that problem \eqref{LP folded concave} is well defined, that is, there
exists a
finite global solution to \eqref{LP folded concave}. To ensure the
well-definedness, it suffices to assume that the statistical loss
function $\mathbb L(\beta)$ is bounded from below on~$\Lambda$.
As we will discuss in Section~\ref{example section}, penalized
linear regression
(least squares), penalized quantile regression, penalized linear
support vector
machine, penalized corrected linear regression and penalized
semiparametric elliptical design regression can all be written in the unified
form of (\ref{LP folded concave}). Thus, the problem setting
in this paper is general enough to cover some new applications that are
not addressed in \citet{Fanetal2012}. Specifically, the
discussions in \citet{Fanetal2012} covered sparse linear
regression, sparse
logistic regression, sparse precision matrix estimation and sparse quantile
regression. All these estimation problems intrinsically have convex loss
functions. \citet{Wangetal2013} and \citet{LohandWainwright}
considered
problems with less regularity by allowing the loss functions to be nonconvex.
Their proposed approaches are, therefore, applicable to corrected linear
regression and semiparametric elliptical design regression.
Nonetheless, both
works assumed different versions of restricted strong convexity. (See
Section~\ref{section nesterov} for more discussions about restricted
strong convexity.)
In contrast,
our analysis does not make assumptions of convexity, nor of any form of
restricted strong convexity, on the statistical loss function.
Moreover, the
penalized support vector machine problem has been addressed in none of the
above literature.

We assume $P_\lambda( \cdot)$ to be either one of the two mainstream
folded concave penalties: (i) smoothly clipped absolute deviation
(SCAD) penalty \citepp{FanandLi2001}, and (ii) minimax concave penalty
[MCP, \citet{Zhang2010}]. Notice that both SCAD and MCP are nonconvex
and nonsmooth.
To facilitate our analysis and computation,
we reformulate \eqref{LP folded concave} into three well-known
mathematical programs: first, a~general quadratic
program; second, a linear program with complementarity constraints; and finally,
a mixed integer (linear) program (MIP). With these reformulations, we
are able to formally state the worst-case complexity
of computing a global optimum to folded concave penalized nonconvex
learning problems.
More importantly, with the MIP reformulation, the global optimal
solution to folded concave penalized nonconvex learning problems can be
numerically solved with a provable guarantee. This reformulation-based solution
technique is referred to as the {\it MIP-based global optimization} (MIPGO).

In this paper,
we make the following major contributions:
\begin{longlist}[(a)]
\item[(a)] We first establish a connection between folded concave
penalized nonconvex learning and quadratic programming. This connection
enables us to analyze the complexity of solving the problem globally.

\item[(b)] We provide an MIPGO scheme (namely, the MIP reformulations)
to SCAD and MCP penalized nonconvex learning, and further prove that
MIPGO ensures global optimality.
\end{longlist}

To our best knowledge, MIPGO probably is the first solution scheme that
theoretically ascertains global optimality. In
terms of both statistical learning and optimization, a global optimization
technique to folded concave penalized nonconvex learning is desirable.
\citet{ZhangandZhang2012} provided a rigorous statement on the statistical
properties of a global solution, while the existing solution techniques in
literature cannot ensure a local minimal solution. Furthermore, the
proposed MIP reformulation enables global optimization
techniques to be applied directly to solving the original nonconvex learning
problem instead of approximating with surrogate subproblems such as local
linear or local quadratic approximations. Therefore, the objective of
the MIP
reformulation also measures the (in-sample) estimation quality. Due to
the critical role of
binary variables in mathematical programming, an MIP has been well
studied in literature.
Although an MIP is theoretically intractable, the computational and algorithmic
advances in the last decade have made an MIP of larger problem scales fairly
efficiently computable \citepp{Bertsimasetal2011}. MIP solvers can further
exploit the advances in computer architectures, for example, algorithm
parallelization, for additional computational power.

To test the proposed solution scheme, we conduct a series of numerical
experiments comparing MIPGO with
different existing approaches in literature. Involved in the comparison
are a local optimization scheme \citepp{LohandWainwright}, approximate
path following algorithm \citepp{Wangetal2013}, LLA \citepp
{Wangetal2013a,Fanetal2012}
and two different versions of coordinate descent algorithms \citepp
{Mazumder2011,BrehenyandHuang2011}.
Our numerical results show that MIPGO can outperform all these
alternative algorithms in terms of solution quality.

The rest of the paper is organized as follows. In Section~\ref{sec2}, we introduce
our setting, present some illustrative examples and derive
reformulations of
nonconvex learning with the SCAD penalty and the MCP in the form of general
quadratic programs. Section~\ref{4} formally states the complexity of
approximating a global optimal solution and then derives MIPGO.
Sections~\ref{section nesterov} and \ref{section LLA}
numerically compare MIPGO with the techniques as per \citet
{Wangetal2013} and \citet{LohandWainwright} and with LLA,
respectively. Section~\ref{statistical performance} presents a more
comprehensive numerical comparison with several existing local schemes.
Section~\ref{conclusion} concludes the paper. Some technical proofs are
given in Section~\ref{some proofs}, and more technical
proofs are given in the online supplement of this paper [\citet{LiuYaoandLi2014}].

%s2 #&#
\section{Setting, example and folded concave penalty reformulation}\label{sec2}

It is worth noting that the abstract form \eqref{LP folded concave}
evidently subsumes a class of nonconvex
learning problems with different statistical loss functions. Before we
pursue further, let us provide a few
examples of the loss functions that satisfy our assumptions to
illustrate the
generality of our statistical setting. Suppose\vspace*{1pt} that $\{(x_t,
y_t):t=1,\dots,n\}\subset\mathbb R^d\times\mathbb R$
is a random sample of size $n$. Let $y=(y_1,\dots,
y_n)^\top$ be the $n\times1$ response vector, and $X=(x_1,\dots,
x_n)^\top$, the $n\times d$ design matrix. Denote throughout this paper
by $\llVert \cdot\rrVert _2$ the $\ell_2$ norm and by $\llvert  \cdot\rrvert  $ the
$\ell_1$ norm.

%s2.1 #&#
\subsection{Examples}\label{example section} (a) The $\ell_2$ loss
for the least squares problem, formulated as
$
\mathbb L_2(\beta):= \frac{1}{2}\sum_{t=1}^n (y_t-x_t^\top\beta
)^2=\frac{1}{2}\llVert  y-X\beta\rrVert ^2_2\nonumber$.
%where $X:= (x_1,\dots, x_n)^\top$.
It is easy to derive that the $\ell_2$-loss can be written in the form
of the loss function as in \eqref{LP folded concave}.

(b) The $\ell_1$ loss, formulated as
$
\mathbb L_1(\beta):= \sum_{t=1}^n \llvert  y_t-x_t^\top\beta\rrvert  =\llvert   y-X\beta
\rrvert  $.
In this case, we can instantiate the abstract form \eqref{LP folded
concave} into
\[
\min_{\beta\in\mathbb R^d,\psi\in\mathbb R^n} \Biggl\{\mathbf{1}^\top \psi+n\sum
_{i=1}^dP_\lambda\bigl(\llvert
\beta_i\rrvert \bigr): -\psi\leq y-X\beta\leq \psi \Biggr\},
\nonumber
\]
where $\mathbf1$ denotes the all-ones vector with a proper dimension.

(c) The quantile loss function in a quantile regression problem,
defined as
\[
\mathbb L_\tau(\beta):= \sum_{t=1}^n
\rho_\tau\bigl(y_t-x_t^\top\beta
\bigr)= \sum_{t=1}^n \bigl(y_t-x_t^\top
\beta\bigr)\bigl\{\tau- \mathbb I\bigl(y_t<x_t^\top
\beta\bigr)\bigr\},
\]
where, for any given $\tau\in(0,1)$, we have $\rho_\tau(u):= u\{\tau
- \mathbb I(u<0)\}$.
This problem with a penalty term can be written in the form of \eqref
{LP folded concave} as
\begin{eqnarray*}
&& \min_{\beta\in\mathbb R^d,\psi\in\mathbb R^n} \mathbf1^\top\bigl[ (y-X\beta ) \tau+\psi
\bigr]+n\sum_{i=1}^dP_\lambda\bigl(
\llvert \beta_i\rrvert \bigr),
\nonumber
\\
&&\qquad\mbox{s.t.}\qquad \psi\geq X\beta-y;\qquad\psi\geq0.\nonumber
\end{eqnarray*}
%
%where $\psi_i$ denotes the $i$-th entry of vector $\psi$.
%{\mathbf Hongcheng: quantile regression can be written as linear program
%as shown at
%http://en.wikipedia.org/wiki/Quantile$ _{-}$regression

%Could you please write the quantile regression?
%}

(d) The hinge loss function of a linear support vector machine
classifier, which is formulated as
$\mathbb L_{SVM}=\sum_{t=1}^n[1-y_tx_t^\top\beta]_+$.
Here, it is further assumed that $y_t\in\{-1,+1\}$, which is the
class label, for all $t=1,\dots,n$. The corresponding instantiation of
the abstract form
\eqref{LP folded concave} in this case can be written as
\begin{eqnarray*}
&& \min_{\beta\in\mathbb R^d,\psi\in\mathbb R^n} \mathbf1^\top\psi +n\sum
_{i=1}^dP_\lambda\bigl(\llvert
\beta_i\rrvert \bigr),
\nonumber
\\
&&\qquad\mbox{s.t.}\qquad\psi_t\geq1-y_tx_t^\top
\beta;\qquad\psi_t\geq0\qquad \forall t=1,\dots,n.
\nonumber
\end{eqnarray*}
%
%{\mathbf Hongcheng: make sure sum is taken over from $i=1$ to $n$ rather
%than $d$.}

(e) Corrected linear regression and semiparametric
elliptical design regression with a nonconvex penalty. According to
\citet{LohandWainwright} and \citet{Wangetal2013} both regression
problems can be
written as general quadratic functions and, therefore, they are
special cases of \eqref{LP folded concave} given that both problems are
well-defined.

%s2.2 #&#
\subsection{Equivalence of nonconvex learning with folded concave
penalty to a general quadratic program}

In this section, we provide equivalent reformulations of the nonconvex
learning problems into a widely investigated
form of mathematical programs, general quadratic programs.
% (Theorem \ref{theorem1} and \ref{reformulate MCP}).
We will concentrate on two commonly-used penalties: the SCAD penalty
and the MCP.
%We will refer to the proposed global optimization method as MIPGO for
%most of the cases.
%Yet, since the proposed reformulation differentiates between solving
%nonconvex learning with the SCAD penalty
%and solving nonconvex learning with the MCP, we will use MIPGO-SCAD or
%MIPGO-MCP to rule out the possible ambiguity occasionally.
%$P_\lambda:\mathbb R_+\rightarrow\mathbb R_+$ is a penalty for
%sparse recovery.
%Note that the two penalties can be defined through non-closed-form
%formulations as follows:

Specifically, given $a>1$ and $\lambda>0$, the SCAD penalty \citepp
{FanandLi2001} is %with $P_{\mathrm{SCAD},\lambda}:\mathbb R_+\rightarrow
%\mathbb R_+$
defined as
%e2.1 #&#
%
\begin{equation}
P_{\mathrm{SCAD},\lambda}(\theta):= \int_{0}^\theta\lambda
\biggl\{\mathbb I(t\leq \lambda)+\frac{(a\lambda-t)_+}{(a-1)\lambda} \mathbb I(t>\lambda) \biggr\}
\,dt,\label{SCAD first order}
\end{equation}
where $\mathbb I( \cdot)$ is the indicator function, and $(b)_+$
denotes the positive part of $b$.
%denotes an index function which takes value 1, if the inequality in
%the argument is satisfied, and which takes value 0, otherwise.

Given $a>0$ and $\lambda>0$, the MCP \citepp{Zhang2010} is %with
%$P_{\mathrm{MCP},\lambda}: \mathbb R_+\rightarrow\mathbb R_+$
defined as
%e2.2 #&#
%
\begin{equation}
P_{\mathrm{MCP},\lambda}(\theta):= \int_{0}^\theta
\frac{(a\lambda-t)_+}{a} \,dt.\label{MCP first order}
\end{equation}
%
%where $( \cdot)_+$ denotes a projection function that yields the
%larger value between 0 and the argument.

We first provide the reformulation of \eqref{LP folded concave} with
SCAD penalty to a general quadratic program
in Proposition~\ref{theorem1}, whose proof will be given in the online
supplement [\citet{LiuYaoandLi2014}]. Let $\mathcal F_{\mathrm{SCAD}}( \cdot, \cdot): \mathbb
R^d\times\mathbb R^d\rightarrow\mathbb R$ be defined as
\[
\mathcal F_{\mathrm{SCAD}}(\beta, g)=\frac{1}{2}\beta^\top Q\beta
\\
+q^\top\beta+n\sum_{i=1}^d
\biggl\{ \bigl(\llvert \beta_i\rrvert -a\lambda\bigr)\cdot
g_i+\frac{1}{2}(a-1)\cdot g_i^2
\biggr\},
\]
where $\beta=(\beta_i)\in\mathbb R^d$ and $g=(g_i)\in\mathbb R^d$.
%%It is easily verifiable that both $\mathcal F_{\mathrm{SCAD}}(\cdot,g)$ and
%$\mathcal F_{\mathrm{SCAD}}(\beta,\cdot)$ are convex.

\begin{proposition}\label{theorem1}
Let $P_\lambda( \cdot)=P_{\mathrm{SCAD},\lambda}( \cdot)$. \textup{(a)} The
minimization problem \eqref{LP folded concave} is equivalent to the
following program:
%e2.3 #&#
%
\begin{equation}
\min_{\beta\in\Lambda, g\in[0,\lambda]^d}\mathcal F_{\mathrm{SCAD}}(\beta, g).\label{biconvex form}
\end{equation}

\textup{(b)} The first derivative of the SCAD penalty can be rewritten as
$P'_\lambda(\theta)=\operatorname{argmin}_{\kappa\in[0, \lambda]}
\{ (\theta-a\lambda)\cdot\kappa+\frac{1}{2}(a-1)\cdot\kappa^2 \}
$ for any $\theta\geq0$.
\end{proposition}

To further simplify the formulation, we next show that program \eqref
{biconvex form}, as an immediate result of
Proposition \ref{theorem1}, is equivalent to a general quadratic program.

\begin{corollary} Program \eqref{biconvex form} is equivalent to
%e2.4 #&#
%e2.5 #&#
%
\begin{eqnarray}
\label{quadratic obj} %\min_{\beta,g,h\in\mathbb R^d} \frac{1}{2}[\beta^\top X^\top X\beta
%+n(a-1)g^\top g+2ng^\top h]+q^\top\beta-na\lambda e^\top g
%\label{quadratic obj}
&& \min_{\beta,g,h\in\mathbb R^d}
\frac{1}{2}\bigl(\beta^\top Q\beta +n(a-1)g^\top
g+2ng^\top h\bigr)+q^\top\beta-na\lambda\mathbf1^\top
g
\nonumber\\[-8pt]\\[-8pt]\nonumber
&&\qquad\mbox{s.t.}\qquad \beta\in\Lambda;\qquad h\geq\beta;\qquad h\geq-\beta;\qquad 0\leq g\leq
\lambda.
\end{eqnarray}
\end{corollary}
\begin{pf}
The proof is completed by invoking Proposition~\ref{theorem1} and the
non-negativity of $g$.
\end{pf}
%
%The program \eqref{QP form 1}-\eqref{QP form 2} can be rewritten in
%the following matrix form

{The above reformulation facilitates our analysis by connecting the
nonconvex learning problem with a general quadratic program. The latter
has been heavily investigated in literature. Interested readers are
referred to \citet{VavasisBook} for an excellent summary on
computational issues in solving a nonconvex quadratic program. }

Following the same argument for the SCAD penalty, we have similar
findings for \eqref{LP folded concave} with the MCP.
The reformulation of \eqref{LP folded concave} with the MCP is given in
the following proposition,
whose proof will be given in the online supplement [\citet{LiuYaoandLi2014}]. Let $\mathcal
F_{\mathrm{MCP}}( \cdot, \cdot): \mathbb R^d\times\mathbb R^d\rightarrow
\mathbb R$ be defined as
\[
\mathcal F_{\mathrm{MCP}}(\beta, g):= \frac{1}{2}\beta^\top Q
\beta+q^\top\beta+n\sum_{i=1}^d
\biggl\{ \frac{1}{2a}g_i^2- \biggl(
\frac{1}{a}g_i-\lambda \biggr)\llvert \beta_i
\rrvert \biggr\}.
\]

{\proposition\label{reformulate MCP} Let $P_\lambda( \cdot
)=P_{\mathrm{MCP},\lambda}( \cdot)$.
\textup{(a)} The model \eqref{LP folded concave} is equivalent to the following program:
%e2.6 #&#
%
\begin{equation}
\min_{\beta\in\Lambda, g\in[0, a\lambda]^d}\mathcal F_{\mathrm{MCP}}(\beta, g).\label{biconvex form MCP}
\end{equation}

\textup{(b)} For any $\theta\geq0$, the first derivative of the MCP can be
rewritten as
$P'_\lambda(\theta)=\frac{a\lambda-g^*(\theta)}{a}$
where $g^*( \cdot):\mathbb R_+\rightarrow\mathbb R$ is defined as
\[
g^*(\theta):= \mathop{\operatorname{argmin}}_{\kappa\in[0, a \lambda
]} \biggl\{
\frac{1}{2a}\kappa^2- \biggl(\frac{1}{a}\kappa-\lambda
\biggr)\theta \biggr\}.
\nonumber
\]
}

Immediately from the above theorem is an equivalence between MCP
penalized nonconvex
learning and the following nonconvex quadratic program.

\begin{corollary} The program \eqref{biconvex form MCP} is equivalent to
%e2.7 #&#
%e2.8 #&#
%
\begin{eqnarray}\label{reformulated LR-MCP obj}
%\min_{\beta,g,h\in\mathbb R^d}~& \frac{1}{2}\beta^\top Q\beta+q^\top
%\beta+n \frac{1}{2a}g^\top g-n\left(\frac{1}{a}g- \lambda\mathbf1
%\right)^\top h
&& \min_{\beta,g,h\in\mathbb R^d}
\frac{1}{2} \biggl(\beta^\top Q\beta+ \frac{n}{a}g^\top
g-\frac{2n}{a}g^\top h \biggr)+ n\lambda\mathbf1^\top
h+q^\top\beta
\nonumber\\[-8pt]\\[-8pt]\nonumber
&&\qquad\mbox{s.t.}\qquad \beta\in\Lambda;\qquad h\geq\beta;\qquad h\geq-\beta;\qquad 0\leq g\leq a
\lambda.%\label{reformulated LR-MCP const}
\end{eqnarray}
\end{corollary}
\begin{pf}
This is a direct result of Proposition~\ref{reformulate MCP} by noting the
non-negativity of $g$.
\end{pf}
With the above reformulations, we are able to provide our complexity
analysis and devise our promised solution scheme.

%s3 #&#
\section{Global optimization techniques}\label{4}
This section is concerned with global optimization of (\ref{LP folded
concave}) with the SCAD penalty and the MCP. We will first establish the
complexity of approximating an $\varepsilon$-suboptimal solution in
Section~\ref{Sec 3} and then provide the promised MIPGO method in
Section~\ref{MIPGO
sec}. Note that, since the proposed reformulation differentiates
between solving nonconvex learning with the SCAD penalty
and solving nonconvex learning with the MCP, we will use MIPGO-SCAD or
MIPGO-MCP to rule out the possible ambiguity occasionally.
%$P_\lambda: \mathbb R_+\rightarrow\mathbb R_+$ is a penalty for
%sparse recovery.
%\subsection{Restricted Strongly Convex Case}
%\citet{LohandWainwright} proposed a linearly convergent method to
%solve
%a regularized M-estimator with nonconvexity, which also applies to
%linear regression with a SCAD/MCP penalty. The proposed method
%requires an assumption of restricted strong convexity together with
%proper choice of parameters. In this section, we show that the
%quadratic program reformulation, under the same condition, is in fact
%strongly convex, and thus, is efficiently solvable and admits linearly
%convergent algorithms.

%\subsection{General Case}
%s3.1 #&#
\subsection{Complexity of globally solving folded concave penalized nonconvex learning}\label{Sec 3}
In Section~\ref{sec2}, we have shown the equivalence between (\ref{LP folded
concave}) and a quadratic program
in both the SCAD and MCP cases. Such equivalence allows us to
immediately apply existing
results for quadratic programs to the complexity analysis of~\eqref{LP
folded concave}. We first
introduce the concept of $\varepsilon$-approximate of global optimum that
will be used in Theorem~\ref{teo31}(c).
Assume that \eqref{LP folded concave} has finite global optimal solutions.
Denote by $\beta^*\in\Lambda$ a finite, globally optimal solution to
\eqref{LP folded concave}.
Following \citet{vavasis}, we call $\beta^*_{\varepsilon}$ to be an
$
\varepsilon$-approximate solution if there exists another
feasible solution $\bar\beta\in\Lambda$ such that
%e3.1 #&#
%
\begin{equation}
\mathcal L\bigl(\beta^*_{\varepsilon}\bigr)-\mathcal L\bigl(\beta^*\bigr)\leq
\varepsilon \bigl[\mathcal L(\bar\beta)-\mathcal L\bigl(\beta^*\bigr)\bigr].
%f(x')-f(x^*)\leq\varepsilon[f(\bar x)-f(x^*)].
\end{equation}

\begin{theorem}\label{teo31}
\textup{(a)} Denote by $I_{d\times d}\in\mathbb R^{d\times d}$ an
identity matrix, and by
%e3.2 #&#
%e3.3 #&#
%e3.4 #&#
%
\begin{equation}
\mathcal H_1:= \lleft[ \matrix{ Q & \mathbf0& \mathbf0
\vspace*{2pt}\cr
\mathbf0& n(a-1)I_{d\times d}&nI_{d\times d}
\vspace*{2pt}\cr
\mathbf0 &nI_{d\times d}&
\mathbf0} \rright]
\end{equation}
the Hessian matrix of \eqref{quadratic obj}. Let $1<a<\infty$, then
$\mathcal H_1$ has at least one negative eigenvalue (i.e., $\mathcal
H_1$
is not positive semidefinite).

\textup{(b)}
Denote by
%e3.5 #&#
%e3.6 #&#
%e3.7 #&#
%
\begin{equation}
\mathcal H_2:= \lleft[ \matrix{ Q & \mathbf0& \mathbf0
\vspace*{2pt}\cr
\mathbf0&
{n}/{a}I_{d\times d}&{-n}/{a}I_{d\times d}
\vspace*{2pt}\cr
\mathbf0 &{-n}/{a}I_{d\times d}&
\mathbf0} \rright]\label{hessian 1}
\end{equation}
the Hessian matrix of \eqref{reformulated LR-MCP obj}. Let $0<a<\infty
$,
then $\mathcal H_2$ has at least one negative eigenvalue.

\textup{(c)} Assume that \eqref{LP folded concave} has finite global
optimal solutions.
Problem \eqref{LP folded concave} admits an algorithm with complexity of
$O(\lceil3d(3d+1)/\sqrt{\varepsilon}\rceil^r l)$ to attain an $\varepsilon
$-approximate of
global optimum, where $l$ denotes the worst-case complexity of
solving a convex quadratic
program with $3d$ variables, and $r$ is the number of negative
eigenvalues of
$\mathcal H_1$ for the SCAD penalty, and the number of negative
eigenvalues of $\mathcal H_2$
for the MCP.
%Denote by $\beta^*\in\Lambda$ a finite,
%globally optimal solution to \eqref{LP folded concave}. %We say that
%$x'$ is an $\varepsilon$-approximate solution
%Then, problem \eqref{LP folded concave} admits an algorithm with
%complexity of
%$O(\lceil3d(3d+1)/\sqrt{\varepsilon}\rceil^r l)$ to find a solution $
%\beta_{\varepsilon}^*\in\Lambda$ such that
%$\mathcal L(\beta^*_{\varepsilon})-\mathcal L(\beta^*)\leq\varepsilon[
%\mathcal L(\bar\beta)-\mathcal L(\beta^*)]$ for some $\bar\beta\in
%\Lambda$. Here $l$ denotes the worst-case complexity of solving a
%convex quadratic
%program with $3d$ variables, and $r$ is the number of negative
%eigenvalues of
%$\mathcal H_1$ for the SCAD case, and the number of negative
%eigenvalues of $\mathcal H_2$
%for the MCP case.
%\item complexity of $O\left(\frac{1}{\varepsilon}\log(\frac{1}{
%\varepsilon}) \log(\log(\frac{1}{\varepsilon}))\right)$ to attain an $
%\varepsilon$-approximation of a KKT point.
\end{theorem}

\begin{pf} Consider an arbitrary symmetric matrix $\Theta$. Throughout
this proof, $\Theta\succeq0$ means that
$\Theta$ is positive semidefinite. % (namely, all eigenvalues of $
%\Theta$ are non-negative).

Notice
$\mathcal H_1 \succeq0$ only if
%e3.8 #&#
%e3.9 #&#
%
\begin{equation}
\mathcal B_1:= \lleft[ \matrix{ n(a-1)I_{d\times d}&nI_{d\times d}
\vspace*{2pt}\cr
nI_{d\times d}&\mathbf0} \rright] \succeq0.
\end{equation}
Since $1<a<\infty$, we have $n(a-1)I_{d\times d}\succ0$. By Schur
complement condition, the positive semidefiniteness of $\mathcal B_1$
requires that
$ -n(a-1)^{-1}\geq0$,
which contradicts with the assumption $1<a<\infty$. Therefore, $
\mathcal H_1$ is not positive semidefinite. This completes the proof
of (a).% Per \citet{PardalosandVavasis} and \citet{Sahni},
%Problem
%\eqref{quadratic obj}-\eqref{quadratic const2} is NP-hard.

%Further noticing that LP-SCAD \eqref{LP-SCAD} can be rewritten as the
%quadratic program \eqref{quadratic obj}-\eqref{quadratic const2}, we
%have the desired result.

In order to show (b), similarly, we have
$\mathcal H_2 \succeq0$ only if
%e3.10 #&#
%e3.11 #&#
%
\begin{equation}
\mathcal B_2:= \lleft[ \matrix{ n/aI_{d\times d}&-n/aI_{d\times d}
\vspace*{2pt}\cr
-n/aI_{d\times d}&\mathbf0} \rright]\succeq0.
\end{equation}
Since $0<a<\infty$, we have $n/aI_{d\times d}\succ0$. By Schur
complement condition, the positive semidefiniteness of $\mathcal B_2$
requires that
$-n/a\geq0$,
which contradicts with the assumption $0<a<\infty$. Therefore, $
\mathcal H_2$ is not positive semidefinite, which means $\mathcal
H_2$ has at least one negative eigenvalue. This completes the proof
for part (b).% Per \citet{PardalosandVavasis} and
%\citet{Sahni}, Problem
%\eqref{quadratic obj}-\eqref{quadratic const2} is NP-hard.

Part (c) can be shown immediately from Theorem 2 in \citet
{vavasis}, and
from the equivalence between \eqref{LP folded concave} and the
quadratic program \eqref{quadratic obj} %-\eqref{quadratic const2}
for the SCAD case, and that between \eqref{LP folded concave} and
\eqref{reformulated LR-MCP obj} %-\eqref{reformulated LR-MCP const}
for the MCP case.
\end{pf}

%{In Theorem~\ref{teo31}(c), solution $\beta_{\varepsilon}^*$ is called an
%$\varepsilon$-approximate of global optimum by \citet{vavasis}.
In Theorem~\ref{teo31}(c), the complexity result for attaining such a
solution is shown in an abstract manner and no practically
implementable algorithm has been proposed to solve a nonconvex
quadratic program in general, or to solve \eqref{quadratic obj} %-
%\eqref{quadratic const2}
or \eqref{reformulated LR-MCP obj} %-\eqref{reformulated LR-MCP const}
in particular.

{\citet{Pardalos} provided an example for a nonconvex quadratic program
with $2^r$ local solutions. Therefore, by the equivalence between
\eqref{quadratic obj} [or \eqref{reformulated LR-MCP obj}] for SCAD (or
MCP) and \eqref{LP folded concave}, the latter may also have $2^r$
local solutions in some bad (not necessarily the worst) cases.}

%s3.2 #&#
\subsection{Mixed integer programming-based global optimization technique}\label{MIPGO sec}
Now we are ready to provide the proposed MIPGO, which essentially is a
reformulation of
nonconvex learning with the SCAD penalty or the MCP into an MIP
problem. Our reformulation is inspired by \citet{MIPliterature}, who
provided MIP reformulations to solve a quadratic program with box constraints.

It is well known that an MIP can be solved with provable global
optimality by solution schemes such as
the branch-and-bound algorithm [B\&B, \citet{marti}].
Essentially, the B\&B algorithm keeps track of both a global lower
bound and a global upper bound on the objective value of the global
minimum. These bounds are updated by B\&B by systematically
partitioning the feasible region into multiple convex subsets and
evaluating the feasible and relaxed solutions within each of the
partitions. B\&B then refines partitions repetitively over iterations.
Theoretically, the global optimal solution is achieved, once the gap
between the two bounds is zero. In practice, the B\&B is terminated
until the two bounds are close enough. The state-of-the-art MIP solvers
incorporate B\&B with additional features such as local optimization
and heuristics to facilitate computation.

%s3.2.1 #&#
\subsubsection{MIPGO for nonconvex learning with the SCAD penalty}

Let us introduce a notation. For two $d$-dimensional vectors $\Phi
=(\phi_i)\in\mathbb R^d$
and $\Delta=(\delta_i)\in\mathbb R^d$, a complementarity constraint $
0\leq\Phi\perp\Delta\geq0$
means that $\phi_i\geq0$, $\delta_i\geq0$, and $\phi_i \delta_i=0$
for all $i:1\leq i\leq d$.
A natural representation of this complementarity constraint is a set of
logical constraints
involving binary variables $\mathbf z=(z_i)\in\{0,1\}^d$:
%\begin{eqnarray}
%\begin{rcases}
%\phi_i\geq0,\qquad\delta_i\geq0,\qquad
%\phi_i \leq\mathcal M z_i
%\\ \delta_i \leq\mathcal M(1- z_i),\qquad z_i\in\{0, 1\}
%\end{rcases}
%\qquad\forall i=1,\dots,d,\label{Mixed integer system reform}
%\end{eqnarray}
%with $\mathcal M$ being a properly large constant. A more succinct
%notation than \eqref{Mixed integer system reform} can be
%e3.12 #&#
%
\begin{eqnarray}\label{Mixed integer system reform succinct}
\Phi\geq0;\qquad
\Delta\geq0;\qquad
\Phi\leq\mathcal M \mathbf z;\qquad
\Delta\leq\mathcal M(1-\mathbf z);\qquad
\mathbf z\in\{0,1\}^d.\hspace*{-30pt}
\end{eqnarray}
The following theorem gives the key reformulation that will lead to the MIP
reformulation.

\begin{theorem}\label{prop 1}
Program \eqref{quadratic obj} %-\eqref{quadratic const2}
is equivalent to a linear program with (linear) complementarity
constraints (LPCC) of the following form:
%e3.13 #&#
%e3.14 #&#
%
\begin{eqnarray}
&& \min \tfrac{1}{2}q^\top\beta-\tfrac{1}{2}\mathbf
b^\top\rho-\tfrac
{1}{2}na\lambda\mathbf1^\top g-
\tfrac{1}{2}\lambda\gamma_4^\top \mathbf1,\qquad\mbox{s.t.}\label{to show comple111}
\\
&&\qquad \cases{Q\beta+q+\gamma_1-\gamma_2+\mathcal A\rho=0;
\qquad ng-\gamma _1-\gamma_2=0,
\vspace*{2pt}\cr
n(a-1)g+nh-na\lambda
\mathbf1-\gamma_3+\gamma_4=0,
\vspace*{2pt}\cr
0\leq
\gamma_{1}\perp h-\beta\geq0; %\label{comple1}
\qquad0\leq
\gamma_{2}\perp h+\beta\geq0,
\vspace*{2pt}\cr
0\leq\gamma_{3}\perp g
\geq0;\qquad0\leq\gamma_{4}\perp\lambda-g\geq0,
\vspace*{2pt}\cr
0\leq\rho\perp
\mathbf b-\mathcal A^\top\beta\geq0,
\vspace*{2pt}\cr
\beta, g, h,
\gamma_1, \gamma_2, \gamma_3,
\gamma_4\in\mathbb R^d;\qquad\rho\in\mathbb
R^m.}\label{comple total 1}
\end{eqnarray}
%
%\textcolor{black}{where for two vectors $u$ and $v$ of the same
%number of
%dimensions, $u\perp v$ means $u_i v_i=0$ for all $i$ with $u_i$ and
%$v_i$ being the $i$-th entry of $u$ and $v$ respectively. }
\end{theorem}

The above LPCC can be immediately rewritten into an MIP. Rewriting the
complementarity constraints in \eqref{comple total 1} % \eqref{comple1}-
%\eqref{comple4}
into the system of logical constraints following \eqref{Mixed integer
system reform succinct}, problem \eqref{quadratic obj} %-
%\eqref{quadratic const2}
now becomes
%
%e3.15 #&#
%e3.16 #&#
%
\begin{eqnarray}
&& \min\tfrac{1}{2}q^\top\beta-\tfrac{1}{2}\mathbf
b^\top\rho-\tfrac
{1}{2}na\lambda\mathbf1^\top g-
\tfrac{1}{2}\lambda\gamma_4^\top \mathbf1;\qquad\mbox{s.t.}\label{mip 1}
\\
&&\qquad \cases{Q\beta+q+\gamma_1-\gamma_2+\mathcal A\rho=0;
\qquad ng-\gamma _1-\gamma_2=0,
\vspace*{2pt}\cr
n(a-1)g+nh-na\lambda
\mathbf1-\gamma_3+\gamma_4=0,
\vspace*{2pt}\cr
\gamma_{1}
\leq\mathcal M \mathbf z_{1};\qquad h_i-
\beta_i\leq\mathcal M(1-\mathbf z_{1}),
\vspace*{2pt}\cr
\gamma_{2}\leq\mathcal M \mathbf z_{2};\qquad h+\beta\leq
\mathcal M (1-\mathbf z_{2}),
\vspace*{2pt}\cr
\gamma_{3}\leq\mathcal M
\mathbf z_{3};\qquad g_i\leq\mathcal M(1-\mathbf
z_{3}),
\vspace*{2pt}\cr
\gamma_{4}\leq\mathcal M \mathbf
z_{4};\qquad\lambda-g\leq\mathcal M (1-\mathbf z_{4}),
\vspace*{2pt}\cr
\rho\leq\mathcal M \mathbf z_5;\qquad\mathbf b-\mathcal
A^\top\beta \leq\mathcal M(1-\mathbf z_5),
\vspace*{2pt}\cr
-\beta\leq
h;\qquad\gamma_2\geq0;\qquad\beta\leq h;\qquad\gamma _1
\geq0,
\vspace*{2pt}\cr
g\geq0;\qquad\gamma_3\geq0;\qquad g\leq\lambda;\qquad
\gamma_4\geq0,
\vspace*{2pt}\cr
\rho\geq0;\qquad\mathbf b-\mathcal
A^\top\beta\geq0,
\vspace*{2pt}\cr
\mathbf z_1, \mathbf z_2,
\mathbf z_3, \mathbf z_4\in\{0, 1\} ^d;
\qquad\mathbf z_5\in\{0, 1\}^m,
\vspace*{2pt}\cr
\beta, g, h,
\gamma_1, \gamma_2, \gamma_3,
\gamma_4\in\mathbb R^d;\qquad\rho\in\mathbb
R^m,} \label{mip 2}
\end{eqnarray}
where we recall that $\mathcal M$ is a properly large constant.

The above program is in the form of an MIP, which admits finite
algorithms that ascertain global optimality.

\begin{theorem}
Program \eqref{mip 1}--\eqref{mip 2} admits algorithms that attain a
global optimum in finite iterations.
\end{theorem}
\begin{pf}
The problem can be solved globally in finite iterations by B\&B \citepp
{LawlerandWood} method.
\end{pf}

The proof in fact provides a class of numerical schemes that solve
\eqref{mip 1}--\eqref{mip 2} globally and finitely. Some of these
schemes have become highly developed and even commercialized. We elect
to solve the above problem using one of the state-of-the-art MIP
solvers, Gurobi, which is a B\&B-based solution tool. (Detailed
information about Gurobi can be found at \surl{http://www.gurobi.com/}.)

%s3.2.2 #&#
\subsubsection{MIPGO for nonconvex learning with the MCP}
Following almost the same argument for the SCAD penalized nonconvex
learning, we can
derive the reformulation of the MCP penalized nonconvex learning
problem into an LPCC per the following theorem. {\theorem
\label{LRmcp}Program \eqref{reformulated LR-MCP obj} is equivalent to the
following LPCC:
%e3.17 #&#
%e3.18 #&#
%
\begin{eqnarray}
\qquad && \min \tfrac{1}{2}q^\top\beta-\tfrac{1}{2}\mathbf
b^\top\rho-\tfrac
{1}{2}a\lambda\mathbf1^\top
\eta_4+\tfrac{1}{2}\lambda n\mathbf1^\top h;\qquad\mbox{s.t.} \label{LPCC 1}
\\
&&\qquad \cases{ Q\beta+q+\eta_1 -\eta_2+\mathcal A\rho=0,
\vspace*{2pt}\cr
\dfrac{n}{a}g-\dfrac{n}{a} h-\eta_3+\eta_4=0;
\qquad-n\biggl(\dfrac{1}{a}g+ \lambda\mathbf1\biggr)-\eta_1-
\eta_2=0,
\vspace*{2pt}\cr
0\leq\eta_1\perp h-\beta\geq0;\qquad0\leq
\eta_2\perp h+\beta\geq0,
\vspace*{2pt}\cr
0\leq\eta_3\perp g\geq0;
\qquad0\leq\eta_4\perp a\lambda\mathbf 1-g\geq0,
\vspace*{2pt}\cr
0\leq\rho\perp
\mathbf b-\mathcal A^\top\beta\geq0,
\vspace*{2pt}\cr
{\beta, g, h,
\eta_1, \eta_2, \eta_3, \eta_4}
\in\mathbb R^d;\qquad\rho\in\mathbb R^m.}\label{LPCC 2}
\end{eqnarray}
}

To further facilitate the computation, program \eqref{LPCC 1}--\eqref
{LPCC 2} can be represented as
%e3.19 #&#
%e3.20 #&#
%
\begin{eqnarray}
\qquad && \min \tfrac{1}{2}q^\top\beta-\tfrac{1}{2}\mathbf
b^\top\rho-\tfrac
{1}{2}a\lambda\mathbf1^\top
\eta_4+\tfrac{1}{2}\lambda n\mathbf1^\top h;\qquad\mbox{s.t.}\label{mip 1 MCP}
\\
&&\qquad \cases{q+Q\beta+\eta_1 -\eta_2+\mathcal A\rho=0;\qquad
\dfrac
{n}{a}g-\dfrac{n}{a} h-\eta_3+\eta_4=0,
\vspace*{2pt}\cr
-n\biggl(\dfrac{1}{a}g+ \lambda\mathbf1\biggr)-\eta_1-
\eta_2=0,
\vspace*{2pt}\cr
0\leq\eta_1\leq\mathcal M \mathbf
z_1;\qquad0\leq h-\beta\leq\mathcal M(1-\mathbf z_1),
\vspace*{2pt}\cr
0\leq\eta_2\leq\mathcal M \mathbf z_2;\qquad0\leq h+\beta
\leq\mathcal M(1-\mathbf z_2),
\vspace*{2pt}\cr
0\leq\eta_3\leq
\mathcal M \mathbf z_3;\qquad0\leq g\leq\mathcal M(1-\mathbf
z_3),
\vspace*{2pt}\cr
0\leq\eta_4\leq\mathcal M \mathbf
z_4;\qquad0\leq a\lambda\mathbf 1-g\leq\mathcal M(1-\mathbf
z_4),
\vspace*{2pt}\cr
0\leq\rho\leq\mathcal M \mathbf z_5;\qquad
\mathbf b-\mathcal A^\top \beta\leq\mathcal M(1-\mathbf
z_5);\qquad\mathbf b-\mathcal A^\top \beta\geq0,
\vspace*{2pt}\cr
\mathbf z_1, \mathbf z_2, \mathbf z_3,
\mathbf z_4\in\{0,1\} ^d;\qquad\mathbf z_5
\in\{0,1\}^m,
\vspace*{2pt}\cr
\eta_1, \eta_2,
\eta_3, \eta_4\in\mathbb R^d;\qquad\rho\in
\mathbb R^m.}\label{mip 2 MCP}
\end{eqnarray}

The computability of a global optimal solution to the above MIP is
guaranteed by the following theorem.

\begin{theorem}
Program \eqref{mip 1 MCP}--\eqref{mip 2 MCP} admits algorithms that
attain a global optimum in finite iterations.
\end{theorem}
\begin{pf}
The problem can be solved globally in finite iterations by B\&B \citepp
{LawlerandWood} method.
\end{pf}

Combining the reformulations in Section~\ref{MIPGO sec}, we want to
remark that the MIP reformulation
connects the SCAD or MCP penalized nonconvex learning with the
state-of-the-art numerical solvers for MIP. This reformulation
guarantees global minimum theoretically and yields reasonable computational
expense in solving \eqref{LP folded concave}. To acquire such a
guarantee, we do not impose very restrictive
conditions. To our knowledge, there is no existing global
optimization technique for the nonconvex learning with the SCAD penalty
or the MCP
penalty in literature under the same or less restrictive assumptions.
More specifically, for MIPGO,
the only requirement on the statistical loss function is that it should
be a lower-bounded quadratic function on the feasible region $\Lambda$
with the Hessian matrix
$Q$ being symmetric. As we have mentioned in Section~\ref{sec2}, an
important class of sparse learning problems
naturally satisfy our assumption. In contrast, LLA, per its equivalence
to a majorization minimization algorithm, converges asymptotically to a
stationary point that does not differentiate among local maxima, local
minima or saddle points. Hence, the resulting solution quality is not
generally guaranteed. \citet{Fanetal2012} proposed the state-of-the-art
LLA variant. It requires restricted eigenvalue conditions to ensure
convergence to an oracle solution in two iterations with a
lower-bounded probability. The convergence of the local optimization
algorithms by \citet{LohandWainwright} and
\citet{Wangetal2013} both require the satisfaction of (conditions that
imply) RSC. To our knowledge, MIPGO stipulates weaker conditions in
contrast to the above solution schemes.

%s3.2.3 #&#
\subsubsection{Numerical stability of MIPGO}\label{section 5}
The representations of SCAD or MCP penalized nonconvex learning
problems as MIPs
introduce dummy variables to the original problem. These dummy
variables are in fact Lagrangian multipliers
in the KKT conditions of \eqref{quadratic obj} or \eqref{reformulated
LR-MCP obj}.
In cases when no finite Lagrangian multipliers exist, the proposed
MIPGO can result in numerical instability.
To address this issue, we study an abstract form of SCAD or MCP
penalized the nonconvex
learning problems given as following:
%e3.21 #&#
%
\begin{equation}
\mathcal\min \bigl\{ \mathcal F(\beta, h, g): {(\beta, h)\in\tilde \Lambda, g
\in[0,M]^d} \bigr\},\label{abstract form1}
\end{equation}
where $M>0$, $\tilde\Lambda:= \{(\beta, h):\beta\in\Lambda,h\in
\mathbb R^d,h\geq\beta,h\geq-\beta\}$, and $\mathcal F:\tilde
\Lambda\times[0,M]^d\rightarrow\mathbb R$ is assumed continuously
differentiable in $(\beta, h, g)$ with the gradient $\nabla
\mathcal F$ being Lipschitz continuous.
It may easily be verified that \eqref{quadratic obj} %-
%\eqref{quadratic const2}
and \eqref{reformulated LR-MCP obj} %-\eqref{reformulated LR-MCP const}
are both special cases of \eqref{abstract form1}.
Now, we can write out the KKT conditions of this abstract problem as:
%e3.22 #&#
%e3.23 #&#
%e3.24 #&#
%e3.25 #&#
%e3.26 #&#
%e3.27 #&#
%
\begin{eqnarray}
\qquad && \nabla_\beta\mathcal F(\beta, h, g)+\upsilon_{1}-\upsilon
_{2}+\mathcal A\rho =0,\label{kkt acs 1}
\\
&& \nabla_h\mathcal F(\beta, h, g)-\upsilon_{1}-
\upsilon_{2} =0,\label{kkt acs 2}
\\
&& \nabla_g\mathcal F(\beta, g)-\zeta_{1}+\zeta_{2}=0,\label{kkt acs 3}
\\
&&\hspace*{-3pt} \lleft. \begin{array} {l}
\zeta_{1,i}\cdot g_i=0;\qquad\zeta_{2,i}
\cdot(g_i-M)=0,
\\
\zeta_{1,i}, \zeta_{2,i}\geq0
\end{array} %
 \rright\}\qquad\forall i=1,\dots,d,\label{comple acs1}
\\
&&\hspace*{-4pt} \lleft. %
\begin{array} {ll} &\upsilon_{1,i}\cdot(
\beta_i-h_i)=0;\qquad\upsilon_{2,i}\cdot(-\beta
_i-h_i)=0,
\\
&\upsilon_{1,i}\geq0;\qquad\upsilon_{2,i}\geq0 \end{array}
 \rright\} \qquad\forall i=1,\dots,d,
\\
&&\rho\geq0,\qquad\rho^\top\bigl(\mathcal A^\top\beta-\mathbf
b\bigr)=0,
\end{eqnarray}
where $\nabla_\beta\mathcal F(\beta, h, g):= {\partial\mathcal
F(\beta, h, g)}/{\partial\beta}$, $\nabla_g\mathcal F(\beta, h,
g):= {\partial\mathcal F(\beta, h, g)}/{\partial g}$, and $\nabla
_h\mathcal F(\beta, h, g):= {\partial\mathcal F(\beta, h,
g)}/{\partial h}$, and where $\zeta_{1},\zeta_{2}, \upsilon
_{1}, \upsilon_{2}\in\mathbb R^d$ and $\rho\in\mathbb R^m$ are the
Lagrangian multipliers that we are concerned with. For convenience, the
$i$th dimension ($i=\{1, \dots,d\}$) of these multipliers are
denoted as $\zeta_{1,i}, \zeta_{2,i}, \upsilon_{1,i}$, $\upsilon
_{2,i}$ and $\rho_i$, respectively. Notice that, since $\mathcal A$
is full-rank, then $\rho$ is bounded if $\llVert \mathcal A\rho\rrVert $ is
bounded, where we let $\llVert \cdot\rrVert $ be an $\ell_p$
norm with
arbitrary $1\leq p\leq\infty$. (To see this, observe that $\llVert
\mathcal A\rho\rrVert _2=\sqrt{\rho^\top\mathcal A^\top\mathcal A\rho}$
and $ \mathcal A^\top\mathcal A$ is positive definite.)

\begin{theorem}\label{stability}
Denote a global optimal solution to problem \eqref{abstract form1} as $
(\beta^*, h^*, g^*)$.
Assume that there exists a positive constant $C_1$ such that
\[
\max\bigl\{\bigl\llVert \nabla_\beta\mathcal F\bigl(
\beta^{*},h^{*},g^{*}\bigr)\bigr\rrVert, \bigl
\llVert \nabla_h\mathcal F\bigl(\beta^{*},h^{*},g^{*}
\bigr)\bigr\rrVert, \bigl\llVert \nabla_g\mathcal F\bigl(
\beta^{*},h^*,g^{*}\bigr)\bigr\rrVert \bigr\}\leq
C_1.
\]
Then the Lagrangian multipliers corresponding to this global optimum,
$\upsilon_1$, $\upsilon_2$, $\zeta_1$, $\zeta_2$ and $\rho$ satisfy that
%e3.28 #&#
%
\begin{equation}
\max\bigl\{\llVert \upsilon_{1}\rrVert, \llVert \upsilon
_{2}\rrVert, \llVert \zeta _{1}\rrVert, \llVert
\zeta_{2}\rrVert \bigr\}\leq C_1; \quad\mbox {and}\quad
\llVert \mathcal A\rho\rrVert \leq3C_1.\label{Assumption Stable}
\end{equation}
\end{theorem}
\begin{pf}
Recall that $\Lambda$ is nonempty. Since $\mathcal A$ is full rank
and all other constraints are linear and nondegenerate, we have the
linear independence constraint qualification satisfied at a global
solution, which then satisfies the KKT condition. (i) In order to show
that $\upsilon_1$ and $\upsilon_2$ are bounded, with \eqref{kkt acs
2}, we have
$\llVert \upsilon_1+\upsilon_2\rrVert = \llVert \nabla_h\mathcal F(\beta
^*,h^*,g^*)\rrVert \leq C_1$.
Noticing the nonnegativity of $\upsilon_1$ and~$\upsilon_2$, we obtain
$\max\{\llVert \upsilon_1\rrVert,\llVert  \upsilon_2\rrVert  \}\leq\llVert \upsilon
_1+\upsilon_2\rrVert = \llVert \nabla_h\mathcal F(\beta^*,h^*,g^*)\rrVert \leq C_1$.
(ii) To show $\mathcal A\rho$ is bounded, considering \eqref{kkt acs 1},
$\llVert \mathcal A \rho\rrVert =\llVert \nabla_\beta\mathcal F
(\beta^*,h^*,g^*)
+\upsilon_1-\upsilon_2\rrVert \leq\llVert \nabla_\beta\mathcal F(\beta
^*,h^*,g^*)\rrVert +\llVert \upsilon_1\rrVert +\llVert  \upsilon_2\rrVert  \leq3C_1$.
(iii) To show that $\zeta_1$ and $\zeta_2$ are bounded, we notice
that, immediately from %\eqref{comple acs}-
\eqref{comple acs1},
$\zeta_{1,i}\geq0; \zeta_{2,i}\geq0$; and $\zeta_{1,i}\cdot\zeta
_{2,i}=0$, for all $i=1,\dots,d$.
Thus, according to \eqref{kkt acs 3},
$C_1\geq\llVert \nabla_g\mathcal F(\beta^*,h^*,g^*)\rrVert \geq\llVert (\max\{
\zeta_{1,i},\zeta_{2,i}\}, i=1,\dots,d)\rrVert $.
Therefore, $\llVert  \zeta_{1}\rrVert  \leq C_1$ and $\llVert \zeta_{2}\rrVert
\leq C_1$.
\end{pf}

{ With Theorem \ref{stability}, we claim that the Lagrangian
multipliers corresponding to a global optimal
solution cannot be arbitrarily large under proper assumptions. Hence,
we conclude that the proposed method
can be numerically stable. In practice, because $\nabla F$ is assumed
Lipschitz continuous, we can simply impose an
additional constraint $\llVert  \beta\rrVert  _\infty\le C$ in the MIP
reformulation for some positive constant
$C$ to ensure the satisfaction of \eqref{Assumption Stable}.
Conceivably, this additional constraint does not result in a
significant modification to the original problem.}
\section{Comparison with the gradient methods}\label{section nesterov}
This section will compare MIPGO with \citet{LohandWainwright} and
\citet{Wangetal2013} when
$\mathbb L:= \mathbb L_{2}$. Thus, the complete formulation is given as
%e4.1 #&#
%
\begin{equation}
\min_{\beta\in\mathbb R^d}\mathcal L(\beta)=\frac{1}{2}\llVert y-X\beta
\rrVert _2^2+n\sum_{i=1}^dP_\lambda
\bigl(\llvert \beta_i\rrvert \bigr),\label{LR specific}
\end{equation}
where $X$ and $y$ are defined as in Section~\ref{example section}.
We will refer to this problem as SCAD (or MCP) penalized linear
regression [LR-SCAD (or -MCP)]. To solve this problem, \citet
{LohandWainwright} and \citet{Wangetal2013} independently
developed two
types of computing procedures based on the
gradient method proposed by
\citet{Nesterov2007}. For the sake of simplicity, we will refer to both
approaches as the gradient methods hereafter,
although they both present substantial differentiation from the
original gradient algorithm proposed by \citet{Nesterov2007}. To
ensure high
computational and statistical performance, both \citet
{LohandWainwright} and
\citet{Wangetal2013} considered conditions called ``restricted
strong convexity''
(RSC). We will illustrate in this section that RSC can
be a fairly strong condition in LR-SCAD or -MCP problems and that MIPGO
may potentially outperform the gradient methods regardless of whether
the RSC is satisfied. %Without such
%conditions, the performance of gradient methods is not guaranteed.
%\textcolor{green}{DELETE: under the assumption of restricted strong
%convexity.
%This section will show by counterexamples that this assumption may not
%be
%satisfied even in some not very irregular cases. In those cases, MIPGO
%can
%outperform both approaches.}

In \citet{LohandWainwright} and \citet{Wangetal2013}, RSC is defined
differently. These two versions of definitions are discussed as below:
%\textbf{\\ Restricted Strong Convexity in
%\citet{LohandWainwright} (RSC
%$_1$)}
let $\beta_{\mathrm{true}}=(\beta_{\mathrm{true},i})\in\mathbb R^d$ be
the true\vspace*{1pt} parameter vector and $k=\llVert \beta_{\mathrm{true}}\rrVert _0$. Denote
that $L(\beta):=\frac{1}{2n}\llVert  y-X\beta\rrVert ^2_2$. Then
according to
\citet{LohandWainwright}, $L(\beta)$ is said to satisfy RSC if the
following inequality holds:
\begin{eqnarray}\label{restricted SC}
\qquad&& L\bigl(\beta'\bigr)- L\bigl(\beta''
\bigr)-\bigl\langle\nabla_\beta L\bigl(\beta''
\bigr),\beta'-\beta ''\bigr\rangle
\nonumber\\[-8pt]\\[-8pt]\nonumber
&&\qquad  \geq \cases{
\displaystyle \alpha_1\bigl\llVert \beta'-
\beta''\bigr\rrVert _2^2-
\tau_1\frac{\log d}{n}\bigl\llvert \beta '-
\beta''\bigr\rrvert ^2, &\quad for all $
\bigl\llVert \beta'-\beta''\bigr\rrVert
_2\leq3$,
\vspace*{2pt}\cr
\displaystyle\alpha_2\bigl\llVert \beta'-
\beta''\bigr\rrVert _2-\tau_2
\frac{\log d}{n}\bigl\llvert \beta '-\beta''
\bigr\rrvert, &\quad for all $\bigl\llVert \beta'-
\beta''\bigr\rrVert _2\geq3$}
\end{eqnarray}
for some $\alpha_1, \alpha_2>0$ and $\tau_1, \tau_2\geq0$.
Furthermore, \citet{LohandWainwright} assumed (in Lemma 3 of their
paper) that $\frac{64k\tau\log d}{n}+\mu\leq\alpha$
with $\alpha=\min\{\alpha_1,\break \alpha_2\}$ and $\tau=\max\{\tau_1,\tau
_2\}$, for some $\mu\geq0$ such that $\mu\llVert  \beta\rrVert  _2^2+\sum_{i=1}^pP_{\lambda}(\llvert  \beta_i\rrvert  )$ is convex. %We show in the
%following that under all these conditions, the nonconvex learning
%problem is strongly convex in some subspaces of the original feasible
%region.

\citet{Wangetal2013} discussed a different version of RSC. They
reformulated \eqref{LR specific} into $\frac{\mathcal L(\beta
)}{n}=\tilde L(\beta)+\lambda\llvert  \beta\rrvert  $, where $\tilde L(\beta
):= L(\beta)+\sum_{i=1}^dP_\lambda(\beta_i)-\lambda\llvert  \beta\rrvert  $.
According to the same paper, one can quickly check that $\tilde L(\beta
)$ is continuously differentiable. Then their version of RSC, as in
Lemma 5.1 of their paper, is given as
%e4.2 #&#
%
\begin{eqnarray}
\tilde L\bigl(\beta'\bigr)-\tilde L\bigl(\beta''
\bigr)\geq\bigl\langle\nabla\tilde L\bigl(\beta ''
\bigr), \beta'-\beta''\bigr\rangle+
\alpha_3\bigl\llVert \beta'-\beta''
\bigr\rrVert ^2_2\label{RSC 2}
\end{eqnarray}
for all $(\beta', \beta'')\in\{(\beta', \beta''):\sum_{i: \beta
_{\mathrm{true},i}= 0}\mathbb I(\beta'_i-\beta''_i\neq0) \leq s\}$ for some $
\alpha_3>0$ and $s\geq k$. Evidently, this implies that \eqref{RSC
2} also holds for all
$\llVert  \beta'-\beta''\rrVert  _0\leq s$.

%of this inequality is implied from the assumptions on the input data
%of the nonconvex learning problem as shown in the same paper.

To differentiate the two RSCs, we will refer to \eqref{restricted SC}
as RSC$_1$, and to
\eqref{RSC 2} as RSC$_2$. A closer observation reveals that both RSCs
imply that the objective function of the nonconvex learning problem is
strongly convex in some sparse subspace involving $k$ number of
dimensions. %More specifically, $\frac{\partial(\mathcal L(\beta)-n
%\lambda\left\vert \beta\right\vert )}{\partial\beta}$ is
%Specifically, if we define $\hat L:\mathbb R^d\rightarrow\mathbb R$
%as $\hat L( \cdot)=L( \cdot)-\mu\left\Vert\cdot\right\Vert_2^2$,
%then the
%following lemma presents their detailed correspondence.

\begin{lemma}\label{RSC 1 to RSC 2}Assume that $L(\beta)$ satisfies
RSC$_1$ in \eqref{restricted SC}.
If $k\geq1$, $\frac{64k\tau\log d}{n}+\mu\leq\alpha$, and $\mu\llVert
\beta\rrVert _2^2+\sum_{i=1}^pP_{\lambda}(\llvert  \beta_i\rrvert  )$ is convex, then
\begin{eqnarray}\label{RSC implied}
&& \frac{1}{n}\mathcal L\bigl(\beta'\bigr)- \frac{1}{n}
\mathcal L\bigl(\beta''\bigr)
\nonumber\\[-8pt]\\[-8pt]\nonumber
&&\qquad \geq \biggl\langle\frac{1}{n}\nabla\mathcal L\bigl(\beta''
\bigr), \beta'-\beta '' \biggr\rangle+
\alpha_3\bigl\llVert \beta'-\beta''
\bigr\rrVert _2^2\qquad  \forall\bigl\llVert
\beta'-\beta''\bigr\rrVert
_0\leq s,
\end{eqnarray}
for some $\alpha_3>0$, where $s=64k-1$,
$
\nabla\mathcal L(\beta'')\in [ n\nabla\tilde L(\beta)+n\lambda
\partial\llvert  \beta\rrvert   ]_{\beta=\beta''}$,
and $\partial\llvert  \beta\rrvert  $ denotes the subdifferential of $\llvert
\beta\rrvert  $.
%\begin{equation}
%\hat L(\beta')-\hat L(\beta'')\geq\nabla\langle\hat L(\beta''),\beta'-
%\beta''\rangle+\tau\log d/n\left\Vert \beta'-\beta''\right\Vert _1^2, \forall\left\Vert
%\beta'-\beta''\right\Vert_0\leq64k-1
%\end{equation}
\end{lemma}

The proof is given in the online supplement [\citet{LiuYaoandLi2014}]. From this lemma, we know that RSC$_1$, together with other
assumptions made by \citet{LohandWainwright}, implies \eqref{RSC
implied} for some $s\geq\llVert \beta_{\mathrm{true}}\rrVert _0=k$ for all $k\geq
1$. Similarly, for RSC$_2$, if the function $\tilde L$ satisfy \eqref
{RSC 2}, in view of the convexity of $\lambda\llvert  \beta\rrvert  $, we
have that
$\frac{\mathcal L(\beta)}{n}=\tilde L(\beta)+\lambda\llvert  \beta\rrvert  $
satisfies \eqref{RSC implied} for some $s\geq\llVert \beta_{\mathrm{true}}\rrVert
_0$. In summary, \eqref{RSC implied} is a necessary condition to both
RSC$_1$ and RSC$_2$.

Nonetheless, \eqref{RSC implied} can be restrictive in some scenarios.
To illustrate this, we conduct a series of simulations as following: we
simulated a sequence of samples $\{(x_t, y_t): 1\leq t\leq n\}$
randomly from the following sparse linear
regression model: $y_t=x^\top_t \beta_{\mathrm{true}}+\varepsilon_t$, for all $
t=1,\dots,n$,
in which $d$ is set to 100, and $\beta_{\mathrm{true}}=[1;1;\mathbf0_{d-2}]$.
Furthermore,
$\varepsilon_t\sim N(0,0.09)$ and $x_t\sim\mathcal N_{d}(0,\Sigma)$
for all $t=1,\dots,n$ with covariance matrix $\Sigma=(\sigma_{ij})\in
\mathbb R^{d\times d}$ defined as
$\sigma_{ij}=\rho^{\llvert  i-j\rrvert  }$. This numerical test considers
only SCAD for an example.
We set the parameters for the SCAD penalty as $a=3.7$ and $\lambda=0.2$.

We conduct a ``random RSC test'' to see if the randomly generated
sample instances can satisfy the RSC condition.
Notice that both versions of RSC dictate that the strong convexity be
satisfied in a sparse subspace that
has only $k$ number of significant parameters. In this example, we have
$k=2$. Therefore, to numerically check if RSC is satisfied,
we conduct the following procedures: (i) we randomly select two
dimensions $i_1, i_2:1\leq i_1<i_2\leq d$;
(ii) we randomly sample two points $\beta^{1}, \beta^2\in\{\beta\in
{\mathbb R}^d:\beta_i=0, \forall i\notin\{ i_1, i_2\}\}$; and
(iii) we check if a necessary condition for \eqref{RSC implied} holds.
That is, we check if the following inequality holds, when $\beta^1\neq
\beta^2$:
%e4.3 #&#
%
\begin{equation}
\frac{\mathcal L(\beta^1)+ \mathcal L(\beta^2)}{2}> \mathcal L \biggl(\frac{\beta^1+\beta^2}{2} \biggr).\label{RSC random check}
\end{equation}
%
%This inequality is a necessary condition for strong convexity in the
%subspace that consists only two significant dimensions.

We consider different sample sizes $n\in\{20, 25, 30, 35\}$ and the
covariance matrix parameters $\rho\in\{0.1, 0.3, 0.5\}$ and
constructed twelve sets of sample instances.
Each set includes 100 random sample instances generated as mentioned
above. For each sample instance,
we conduct 10,000 repetitions of the ``random RSC test.'' If \eqref
{RSC random check} is satisfied for all these 10,000 repetitions,
we say that the sample instance has passed the ``random RSC test.''
Table~\ref{random check 1a} reports the test results.

\begin{figure}[b]

\includegraphics{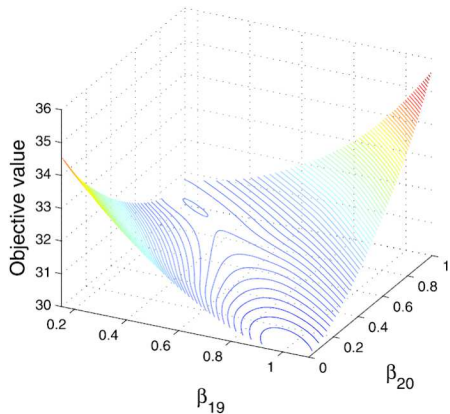}

\caption{A sample instance that fails in the random RSC test.}
\label{RSC test typical instance}
\end{figure}

We observe from Table~\ref{random check 1a} that in some cases the
percentage for passing the random RSC
test is noticeably low. However, with the increase of sample size, that
percentage grows quickly. Moreover,
we can also observe that when $\rho$ is larger, it tends to be more
difficult for RSC to hold.
Figure~\ref{RSC test typical instance} presents a typical instance that
does not satisfy RSC
when $n=20$ and $\rho=0.5$. This figure shows the 3-D contour plot of
objective function when
the decision variable is within the subspace $\{\beta: \beta
_i=0, \forall i\notin\{19, 20\}\}$.
We can see that the contour plot apparently indicates nonconvexity of
the function in the subspace,
which violates \eqref{RSC random check}.

%t1 #&#
%
\begin{table}[t]
\tabcolsep=0pt
\tablewidth=250pt
\caption{Percentage for successfully passing the random RSC test out of
100 randomly generated instances}\label{random check 1a}
\begin{tabular*}{\tablewidth}{@{\extracolsep{\fill}}@{}lcccc@{}}
\hline
$\bolds{\rho}$ &$\bolds{n=35}$ &$\bolds{n=30}$ &$\bolds{n=25}$ &$\bolds{n=20}$\\
\hline
0.1 &93\% &81\% &53\% &4\% \\
0.3 &94\% &76\% &39\% &9\% \\
0.5 &55\% &50\% &21\% &1\% \\
\hline
\end{tabular*}
\end{table}

%f1 #&#
%

%\begin{minipage}{0.2\textwidth}
%\begin{tabular}{\right\vert c\right\vert c\right\vert c\right\vert }
%\hline
% A & B & C \\
%\hline
% 1 & 2 & 3 \\
%\hline
% 4 & 5 & 6 \\
%\hline
%\end{tabular}
%\end{minipage}
%\begin{minipage}{0.2\textwidth}
%\begin{tabular}{c\right\vert c\right\vert c}
% A & B & C \\
%\hline
% 1 & 2 & 3 \\
%\hline
% 4 & 5 & 6 \\
%\end{tabular}
%\end{minipage}

%t2 #&#
%
\begin{table}%[h!]
\tabcolsep=0pt
\caption{Comparison between MIPGO and the gradient methods.
%``AD" denotes the average absolute deviation; ``FP" denotes the
%average number of ``false positive";
%``FN" denotes the average number of ``false negative"; ``Gap'' denotes
%the absolute difference from the objective
%function of the MIPGO solution; and, lastly, ``Time'' stands for the
%computational time.
The numbers in parenthesis are the standard errors. GM$_1$ and GM$_2$
stand for the gradient methods proposed
by \citet{LohandWainwright} and \citet{Wangetal2013}, respectively}\label{table 3new}
\begin{tabular*}{\tablewidth}{@{\extracolsep{\fill}}lccccc@{}}\hline
\textbf{Method}  & \textbf{AD} & \textbf{FP} &\textbf{FN} &\textbf{Gap} & \textbf{Time}\\
\hline
\multicolumn{6}{@{}c@{}}{$\rho=0.5$, $n=20$}\\
MIPGO & 0.188 & 0.230 & 0 & 0 & 29.046 \\
&(0.016) & (0.042) & (0) & (0) & (5.216) \\[3pt]
GM$_1$& 2.000 & 0 & 2&25.828 &0.002 \\
&(0.000) & (0) & (0) & (0.989)&(0.001) \\[3pt]
{GM$_2$}& 0.847 & 5.970 & 0 & 1.542 & 0.504 \\
&(0.055) & (0.436) & (0) & (0.119) & (0.042) \\[6pt]
\multicolumn{6}{@{}c@{}}{$\rho=0.1$, $n=35$}\\
{MIPGO}& 0.085 & 0.020 & 0 & 0 &27.029 \\
&(0.005) & (0.141) & (0) & (0) &(4.673) \\[3pt]
{GM$_1$}& 2.000 & 0 & 2 &31.288 & 0.002 \\
&(0.000) & (0) & (0) &(1.011)& (0.000) \\[3pt]
{GM$_2$}& 0.936 & 6.000 & 0 & 4.179& 0.524 \\
&(0.044) & (0.348) & (0) & (0.170) & (0.020) \\
\hline
\end{tabular*}
\end{table}

We then compare MIPGO with both gradient methods in two sets of the
sample instances from the table:
(i) the one that seems to provide the most advantageous problem
properties ($\rho=0.1$, and $n=35$) to
the gradient methods; and (ii) the one with probably the most
adversarial parameters ($\rho=0.5$, and $n=20$) to the gradient
methods. Notice that the two gradient methods are implemented on MatLab
following the descriptions by \citet{Wangetal2013} and \citet
{LohandWainwright}, respectively, including their initialization
procedures. MIPGO is also implemented on MatLab calling Gurobi
(\surl{http://www.gurobi.com/}). We use CVX, ``a package for specifying and
solving convex programs'' [\citeauthor{CVX1} (\citeyear{CVX1,CVX2})], as the interface between
MatLab and Gurobi.
Table~\ref{table 3new} presents our comparison results in terms of
computational, statistical and optimization measures.
More specifically, we use the following criteria for our comparison:
\begin{itemize}
\item Absolute deviation (AD), defined as the distance between the
computed solution and the true parameter vector. Such a distance is
measured by $\ell_1$ norm.
\item False positive (FP), defined as the number of entries in the
computed solution that are wrongly selected as nonzero dimensions.
\item False negative (FN), defined as the number of entries in the
computed solution that are wrongly selected as zero dimensions.
\item Objective gap (``Gap''), defined as the difference between the
objective value of the computed solution and the objective value of the
MIPGO solution. A positive value indicates a worse relative performance
compared to MIPGO.
\item Computational time (``Time''), which measures the total
computational time to generate the solution.
\end{itemize}
AD, FP and FN are commonly used statistical criteria, and ``Gap'' is a
natural measure of optimization performance.
In Table~\ref{table 3new}, we report the average values for all the
above criteria out of 100 randomly generated instances aforementioned.
From this table, we observe an outperformance of MIPGO over the other
solution schemes on solution quality
for both statistical and optimization criteria. However, MIPGO
generates a higher computational
overhead than the gradient methods.

\section{Numerical comparison on optimization performance with local
linear approximation}\label{section LLA}
In this section, we numerically compare MIPGO with local
linear approximation (LLA). We implement LLA on MatLab. In the
implementation, we invoke the procedures of LLA iteratively until the
algorithm fully converges. This shares the same spirit
as the multistage procedure advocated by \citet
{HuangandZhang2012}. At
each iteration, the LASSO subproblem is solved with Gurobi 6.0 using
CVX [\citeauthor{CVX1} (\citeyear{CVX1,CVX2})] as the interface.
We report in the following a series of comparison results
in terms of the optimization accuracy.

%f2 #&#
%
\begin{figure}

\includegraphics{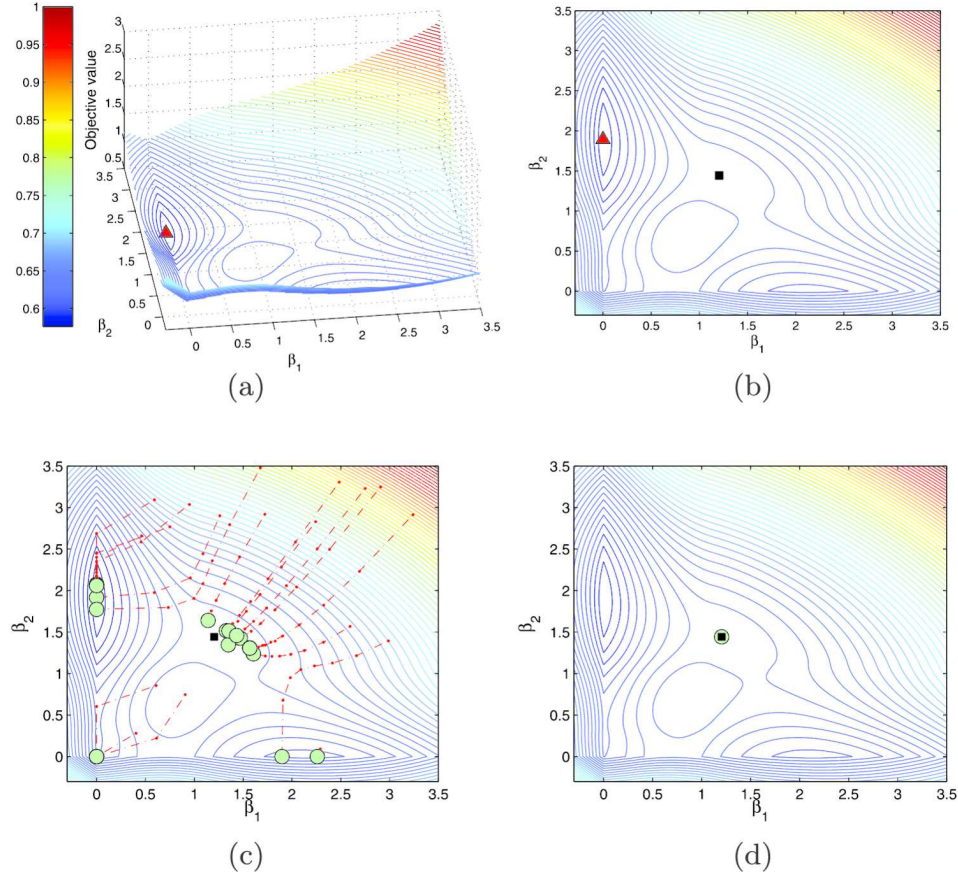}

\caption{\textup{(a)} 3-D contour plots of the 2-dimension LR-SCAD problem and
the solution generated by MIPGO in
20 runs with random initial solutions. The triangle is the MIPGO
solution in both
subplots. \textup{(b)} 2-D representation of subplot \textup{(a)}. \textup{(c)} Trajectories of 20
runs of LLA with random initial solutions.
\textup{(d)} Trajectories of 20 runs of LLA with the least squares solution as
the initial solutions.}\label{SCAD small}
\end{figure}

%f3 #&#
%
\begin{figure}

\includegraphics{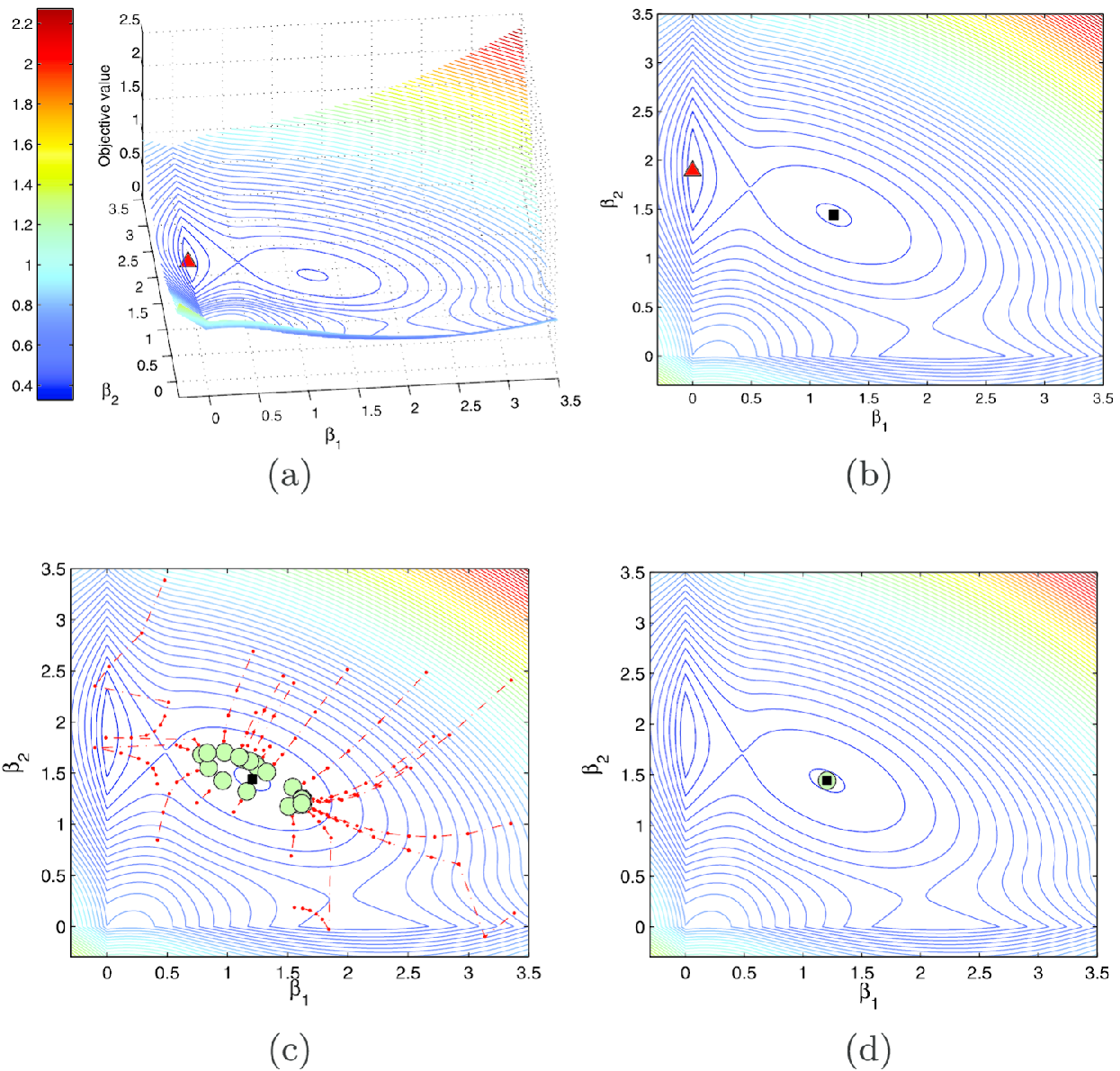}

\caption{\textup{(a)} 3-D contour plots of the 2-dimension LR-MCP problem and
the solution generated by MIPGO in 20 runs with random initial
solutions. The triangle is the MIPGO solution in both subplots. \textup{(b)}
2-D representation of subplot \textup{(a)}. \textup{(c)} Trajectories of 20 runs of LLA
with random initial solutions. \textup{(d)} Trajectories of 20 runs of LLA with
the least squares solution as the initial solutions.}\vspace*{-10pt}\label{MCP small}
\end{figure}

%\subsection{Numerical Test on Optimization Criterion}
%s5.1 #&#
\subsection{Numerical tests on a two-dimensional problem}\label
{2-dimension LLA}
In the following, we conduct a numerical test on a two-dimensional
LR-SCAD and
a two-dimensional LR-MCP problem. We generate one instance for both of
LR-SCAD and LR-MCP problems through the following procedures: we
randomly generate $\beta_{\mathrm{true}}\in\mathbb R^2$ with a uniformly
distributed random vector on $[-1, 5]^2$ and
then generate 2 observations $x_t\sim\mathcal N_2(0,\Sigma)$, $t\in\{
1, 2\}$, with
covariance matrix $\Sigma=(\sigma_{ij})$ and $\sigma_{ij}=0.5^{\llvert  i-j\rrvert  }$. Finally, we compute $y_t$ as $y_t= x_t^\top\beta+
\varepsilon_t$
with $\varepsilon_t\sim\mathcal N(0,1)$ for all $t\in\{1,2\}$. Both the
LR-SCAD problem and the LR-MCP problem use the same set of samples $\{
(x_t,y_t): t=1, 2\}$ in their statistical loss functions. The only
difference between the two is the different choices of penalty
functions. The parameters for the
penalties are prescribed as $\lambda=1$ and $a=3.7$ for the SCAD and
$\lambda=0.5$ and $a=2$ for the MCP. Despite their small dimensionality,
these problems are nonconvex with multiple local solutions. Their nonconvexity
can be visualized via the 2-D and 3-D contour plots provided in
Figure~\ref{SCAD small}(a)--(b) (LR-SCAD) and Figure~\ref{MCP
small}(a)--(b) (LR-MCP).

%, and (iv) initial value set to be the LASSO solution (denoted by
%LLA$_1$).
We realize that the solution generated by LLA may depend on its
starting point. Therefore, to make a fair numerical comparison, we
consider two possible initialization procedures: (i) LLA with random
initial solutions generated with a uniform distribution on the set $[0,
3.5]^2$ (denoted LLA$_r$), and (ii) LLA with initial solution set to
be the least squares solution (denoted by LLA$_{\mathrm{LSS}}$). (We will also
consider LLA initialized with LASSO in later sections.)

To fully study the impact of initialization to the solution quality, we
repeat each solution scheme 20 times. The best (Min.), average (Ave.)
and worst (Max.)
objective values as well as the relative
objective difference (gap$_{(\%)}$) obtained in the 20 runs are
reported in Table~\ref{table1}.
Here, gap$_{(\%)}$ is defined
\begin{eqnarray}
\frac{\{\mbox{Objective of computed solution}\}-\{\mbox{Objective of
MIPGO solution}\}}{\{\mbox{Objective of computed solution}\}}\times 100\%.
\nonumber
\end{eqnarray}
%
% In specific, the relative objective difference is calculated as in
%\eqref{Gap}.
%
%where $\beta_{LLA_x}$ denotes the solution generated by an LLA
%approach with different initialization schemes, $\beta_{MIP}$
%denotes the solution generated by MIPGO approach, and $\mathcal{L}(
%\cdot)$ is as defined in \eqref{LP folded concave}. Notice that a
%positive value of $gap_{(\%)}$ indicates that the solution generated
%by MIPGO has a better objective value than that by LLA approach with
%the corresponding initialization scheme.
From the table, we have the following observations:
\begin{longlist}[4.]
\item[1.] LLA$_r$'s performance varies in different runs. In the best
scenario, LLA attains the global optimum, while the average performance
is not guaranteed.
\item[2.] LLA$_{\mathrm{LSS}}$ fails to attain the global optimal solution.
\item[3.] LLA with either initialization procedure yields a local
optimal solution.
\item[4.] MIPGO performs robustly and attains the global solution at
each repetition.
\end{longlist}

%t3 #&#
%
\begin{table}
\tabcolsep=0pt
\caption{Test result on a toy problem.
``gap$_{(\%)}$'' stands for the relative difference in contrast to MIPGO} \label{table1}
\begin{tabular*}{\tablewidth}{@{\extracolsep{\fill}}@{}lcccccc@{}}
\hline
\textbf{Penalty} & \textbf{Measure} &\textbf{LLA}$_{\bolds{r}}$ & \textbf{gap}$_{\bolds{(\%)}}$ &\textbf{LLA}$_{\mathbf{LSS}}$& \textbf{gap}$_{\bolds{(\%)}}$&
\textbf{MIPGO}\\
\hline
SCAD
&Min & 0.539 &\phantom{0}0.00&0.900&40.12& 0.539\\
&Ave &0.911&40.85&0.900&40.12& 0.539\\
&Max &2.150&74.93&0.900& 40.12& 0.539\\[3pt]
MCP
&Min & 0.304 &\phantom{0}2.63&0.360&17.78& 0.296\\
&Ave &0.435&31.95&0.360&17.78& 0.296\\
&Max &1.293&77.11&0.360& 17.78& 0.296\\\hline
\end{tabular*}
\end{table}

Figures~\ref{SCAD small}(c) and \ref{MCP small}(c) present the
search trajectories
(dot dash lines) and convergent points (circles) of LLA$_r$ for
LR-SCAD and LR-MCP,
respectively. In both figures, we observe a high dependency of LLA's
performance on the
initial solutions. Note that the least squares solutions for the two
problems are denoted by
the black squares. Figures~\ref{SCAD small}(d) and \ref{MCP
small}(d) present the convergent
points of LLA$_{\mathrm{LSS}}$ for \mbox{LR-SCAD} and for LR-MCP, respectively. LLA$
_{\mathrm{LSS}}$ utilizes the least squares
solution (denoted as the black square in the figure) as its starting
point. This least squares solution
happens to be in the neighborhood of a local solution in solving both
problems. Therefore,
the convergent points out of the 20 repetitions of LLA$_{\mathrm{LSS}}$ all
coincide with the
least squares solution. Even though we have $n=d=2$ in this special
case, we can see
that choosing the least squares solution as the initial solution may
lead the LLA to a
nonglobal stationary point. The solutions obtained by MIPGO is
visualized in Figure~\ref{SCAD small}(b)
and \ref{MCP small}(b) as triangles. MIPGO generates the same solution
over the 20 repetitions even with random initial points.
%({\mathbf Hongcheng: What are the value of $\lambda$'s used in SCAD and
%MCP})

%s5.2 #&#
\subsection{Numerical tests on larger problems}
In the following, we conduct similar but larger-scale simulations to
compare MIPGO and LLA in terms of optimization performance.
For these simulations, we randomly generate problem instances as
follows: we first
randomly generate a matrix $T\in\mathbb R^{d\times d}$ with the entry
on $i$th row and
$j$th column uniformly distributed on $[0,  0.5^{\llvert  i-j\rrvert  }]$ and
set $\Sigma=T^\top T$ as the covariance matrix. % with $
%\sigma_{ij}=0.5^{\left\vert i-j\right\vert }$;
Let the true parameter vector $\beta
_{\mathrm{true}}=[3~2~10~0~1~1~2~3~1.6~6~\mathbf0_{1\times(d-10)}]$. We then
randomly generate a sequence of observations $\{(x_t, y_t): t=1,\dots,n\}$ following a linear model $y_t=x_t^\top\beta_{\mathrm{true}} +\varepsilon_t$, where
$x_t \sim\mathcal N_d(\mathbf0, \Sigma)$, and $\varepsilon_t\sim
\mathcal N(0,1.44)$ for all $t=1,\dots, n$. Finally, the penalty
parameters are $\lambda=1$ and $a=3.7$ for SCAD and $\lambda=0.5$
and $a=2$ for MCP.

%We also consider LLA with different initialization schemes. More
%specifically, we involve in our comparison both LLA$_r$ as in
%Subsection \ref{2-dimension LLA} and the LLA taking the final solution
%of LASSO as the initial solution (denoted by LLA$_1$). Notice that
%LLA$_{LLS}$ is a special case of LLA$_1$ when the coefficients of the
%LASSO penalty is zero.

Following the aforementioned descriptions, we generate problem
instances with different problem sizes $d$ and sample sizes $n$ (with
3 problem instances generated for each combination of $d$ and $n$) and
repeat each solution scheme 20 times. For these 20 runs, we randomly
generate initial solutions for MIPGO with each entry following a
uniform distribution on $[-10, 10]$. Similar to the 2-dimensional
problems, we also involve in the comparison LLA with different
initialization procedures:
(i) LLA with randomly generated initialization solution whose each
entry follows a uniform distribution on $[-10, 10]$ (denoted LLA$_r$).
(ii) LLA with zero vector as the initial solution (denoted LLA$_0$).
(iii) LLA with the initial solution prescribed as the solution to the
LASSO problem (denoted LLA$_1$). More specifically, the LASSO problem
used in the initialization of LLA$_1$ is formulated as
%e5.1 #&#
%
\begin{equation}
\min_{\beta\in\mathbb R^n}\frac{1}{2}\llVert y-X\beta\rrVert
_2^2+\omega\sum_{i=1}^d
\llvert \beta_i\rrvert, \label{LASSO problem}
\end{equation}
where $X:= (x_1,\dots,x_n)^\top$, $y:= (y_1,\dots,y_n)^\top$, and $
\omega:= \frac{1}{10}n\lambda\cdot(K-1)$ at $K$th run with $
1\leq
K\leq20$.
This is designed to examine how sensitive the performance the LLA$_1$
depends on the initial estimate.
We would also like to remark that, when $K=1$, the initial solution
for LLA will be exactly the least squares solution. We would also like
to remark that the LLA initialized with LASSO is the solution scheme
proposed by \citet{Fanetal2012}.

%t4 #&#
%
\begin{table}
\tabcolsep=0pt
\caption{Numerical comparison of LLA and the proposed MIPGO on LR-SCAD and LR-MCP problems
with different problem scales. ``TS'' stands for ``Typical Sample''}\label{table 2}
\begin{tabular*}{\tablewidth}{@{\extracolsep{\fill}}@{}lcd{3.2}d{3.3}d{2.2}d{3.2}d{2.2}d{4.2}d{2.2}@{}}
\hline
& & \multicolumn{1}{c}{\textbf{MIPGO}} & \multicolumn{1}{c}{\textbf{LLA}$_{\bolds{r}}$}
& \multicolumn{1}{c}{\textbf{gap}$_{\bolds{(\%)}}$} & \multicolumn{1}{c}{\textbf{LLA}$_{\bolds{0}}$}
& \multicolumn{1}{c}{\textbf{gap}$_{\bolds{(\%)}}$} & \multicolumn{1}{c}{\textbf{LLA}$_{\bolds{1}}$}
& \multicolumn{1}{c@{}}{\textbf{gap}$_{\bolds{(\%)}}$} \\
\hline
\multicolumn{9}{@{}c@{}}{LR-SCAD}\\
TS 3   & Min. & 89.87 & 89.87 &  0.00 &104.96 & 14.37 & 104.96 & 14.37 \\
$d=10$& Ave. & 89.87 & 109.19 & 17.69 & 104.96 & 14.37 & 104.96 &14.37 \\
$n=10$& Max. & 89.87 & 162.93 & 44.84 & 104.96 & 14.37 & 104.96 &14.37
\\[3pt]
TS 6  & Min. &86.04 & 88.219 & 2.46 &115.17 &25.30 &108.37 &20.60\\
$d=20$ & Ave. & 86.04 & 105.13 & 18.15 &115.17&25.30 &108.37 &20.60\\
$n=10$& Max. &86.04 & 143.51 & 40.05 & 115.17&25.30 &108.37 &20.60
\\[3pt]
TS 9   & Min. & 120.35 & 120.35 &  0.00 & 150.15&19.85 & 120.35 & 0.00\\
$d=40$ & Ave. & 120.35 & 167.21 & 28.02 & 150.15&19.85& 120.35 & 0.00\\
$n=15$ & Max. & 120.35 & 203.18 & 40.76 & 150.15&19.85& 120.35 & 0.00
\\[3pt]
TS 12   & Min. &519.14 & 519.14 &  0.00 & 560.28 &  7.34  &538.47 &   3.59 \\
$d=200 $  & Ave. &519.14 & 733.06 &  29.18  & 560.28 & 7.34 & 538.47 &   3.59 \\
$n=60$ & Max. &519.14 & 959.00 &  45.87  & 560.28 & 7.34 & 538.47 &   3.59
\\[3pt]
TS 15   & Min. & 841.72 & 841.90 & 0.02 &1003.69 & 16.14 &873.44 &3.63\\
$d=500$  & Ave. &841.72 & 981.73 & 14.26 &1003.69 & 16.14 &873.44 &3.63\\
$n=80$ & Max. &841.72 & 1173.06 & 28.25 &1003.69 & 16.14 &873.44 & 3.63
\\[3pt]
TS 18   & Min. & 1045.22 & 1105.70 &  5.47  &1119.84 & 6.66 &1119.84&6.66\\
$d=1000$  & Ave. &1045.22 & 1135.84 & 7.98 &1119.84 & 6.66 &1119.84 &6.66\\
$n=100$ & Max. & 1045.22 & 1309.70& 20.19 &1119.84 & 6.66 &1119.84 &6.66
\\[3pt]
\multicolumn{9}{@{}c@{}}{LR-MCP}\\
TS 3   & Min. & 13.65 & 15.77 &13.43 &21.51 & 36.54 & 25.00 & 45.39 \\
$d=10 $ & Ave. & 13.65 & 20.60 & 39.59 & 21.51 & 36.54 & 25.00&45.39 \\
$ n=10$ & Max. & 13.65 & 32.83 & 58.41 & 21.51 & 36.54 & 25.00 &45.39
\\[3pt]
TS 6  & Min. &14.71 & 17.60 & 16.41 &14.71 & 0.00 &20.06 &26.67\\
$d=20 $  & Ave. & 14.71 & 17.60 & 22.54 &14.71 & 0.00 &20.06 &26.67\\
$ n=10$ & Max. &14.71 & 17.60 & 65.38 &14.71 & 0.00 &20.06 &26.67
\\[3pt]
TS 9 & Min. &23.64 & 27.17 & 13.02 & 26.57&11.05 & 49.08 &51.84\\
$d=40 $  & Ave. &23.64 & 27.17 & 35.40 & 26.57&11.05 & 49.08 &51.84\\
$ n=15$ & Max. &23.64 & 27.17 & 57.98 & 26.57&11.05 & 49.08 &51.84
\\[3pt]
TS 12  & Min. & 93.55 & 105.62&  11.42  &112.13 & 16.57 &120.63 & 22.45\\
$d=200 $  & Ave. & 93.55 & 165.25 &  43.39  & 112.13 & 16.57 & 120.63 & 22.45\\
$ n=60$ & Max. & 93.55 & 596.52 &  84.32  &112.13 &16.57 & 120.63 & 22.45
\\[3pt]
TS 15  & Min. & 163.98 & 175.44 & 6.53 &221.84 & 26.08 &179.53 &8.66\\
$d=500 $  & Ave. &163.98 & 211.62 & 22.51 &221.84 & 26.08 &179.53&8.66\\
$ n=80$ & Max. &163.98 & 237.56 & 30.97 &221.84 & 26.08 &179.53 &8.66
\\[3pt]
TS 18  & Min. & 249.89 & 267.83 &  6.70  &267.83 & 6.70 & 272.39&8.27\\
$d=1000 $  & Ave. & 249.89 & 322.24 & 22.25 &267.83 & 6.70 & 272.39 &8.27\\
$ n=100$& Max. & 249.89 & 530.60& 52.60 &267.83 & 6.70 & 272.39 &8.27
\\
\hline
\end{tabular*}
\end{table}

The best (Min.), the average (Ave.) and the worst (Max.) objective
values and the
relative objective differences (gap$_{(\%)}$) of the 20 runs for each
instance are reported in
the upper and lower panels of Table~\ref{table 2} for LR-SCAD and
LR-MCP, respectively.
Notice that for
each problem scale, we generate three test instances randomly, but
Table~\ref{table 2} only reports one of the three
instances for each problem size due to the limit of space.
Tables S1 %\ref{table 2 complement}
and S2 %\ref{table 2 MCP complement}
in Appendix S4 %\ref{Additional test tables}
will complement the rest of the results.
According to the numerical results, in all instances with different
dimensions, MIPGO yields the lowest objective value, and in many cases,
gap$_{(\%)}$ value is nontrivially large. This indicates the
outperformance of our proposed MIPGO over all counterpart algorithms.
%\item In many cases, the Gap.
%\item MIPGO outperforms LLA in every problem instance with significant
%improvement on the average and the worst objective values.
%2) LLA with random initializations (LLA$_{r}$) attains zero gap(\%)
%value in the best cases in solving most instances. Nonetheless, the
%average and the worst performance are significantly worse than the
%MIPGO. This indicates that, contingent on a proper initial solution,
%the LLA converges to the global solution. %This verifies the
%theoretical findings presented collectively through Theorem \ref{DAM
%conv} and \ref{LLA conv} that the LLA approach attains a KKT point.
%3) For larger dimensional problems, (when $d\geq40$), the LLA with
%LASSO solutions as initial solutions has smaller gap$_{(\%)}$ values
%than the LLA with zero vector as initial solution, indicating a
%tendency of LLA$_1$'s improvement over LLA$_0$ on solution quality
%in problems with larger dimensions.

%s6 #&#
\section{Numerical comparison on statistical performance with local algorithms}\label{statistical performance}
We next examine MIPGO on the statistical performance in comparison with
several existing local algorithms, including coordinate descent, LLA,
and gradient methods.
We simulate the random samples $\{(x_t, y_t), t=1,\dots,n\}$ from
the following linear
model$:y_t=x_t^\top(\beta_{\mathrm{true},i}: 1\leq i\leq d-1)+\beta
_{\mathrm{true},d}+\varepsilon_t$,
where we let $d=1001$, $n=100$, and $\beta_{\mathrm{true},d}$ is the intercept.
$\beta_{\mathrm{true}}$ is constructed by first setting $\beta_{\mathrm{true},d}=0$,
then randomly choosing 5 elements among dimensions $\{1,\dots,d-1\}$ to
be 1.5, and setting the other $d-6$ elements as zeros. Furthermore, for
all $t=1,\dots, n$, we let
$\varepsilon_t\sim N(0,1.44)$ and $x_t\sim\mathcal N_{d-1}(0,\Sigma)$
with $\Sigma=(\sigma_{ij})$ defined as
$\sigma_{ij}=0.5^{\llvert  i-j\rrvert  }$. For both LR-SCAD and LR-MCP, we
set the parameter $a=2$, and tune $\lambda$ the same way as presented
by \citet{Fanetal2012}. We generate 100 instances using the above
procedures, and solve each of these instances using MIPGO and other
solutions schemes, including: (i) coordinate descent; (ii)~gradient
methods; (iii) SCAD-based and MCP-based LLA; and (iv) the LASSO method.
The relative details of these techniques are summarized as follows:
\begin{longlist}
\item[\textit{LASSO}:] The LASSO penalized linear regression, coded in
MatLab that invokes Gurobi 6.0 using CVX as the interface.
\item[\textit{GM}$_1$-\textit{SCAD}/\textit{MCP}:] The SCAD/MCP penalized linear regression
computed by the local solution method by \citet{LohandWainwright}
on MatLab.
\item[\textit{GM}$_2$-\textit{SCAD}/\textit{MCP}:] The SCAD/MCP penalized linear regression
computed by the approximate path following algorithm by \citet
{Wangetal2013} on MatLab.
\item[\textit{SparseNet}:] The R-package {\it sparsenet} for SCAD/MCP
penalized linear regression computed by coordinate descent \citepp{Mazumder2011}.
\item[{\it Ncvreg-SCAD}/\textit{-MCP}:] The R-package {\it ncvreg} for MCP
penalized linear regression computed by coordinate descent \citepp
{BrehenyandHuang2011}.
\item[{\it SCAD-LLA$_1$}/\textit{MCP-LLA$_1$}:] The SCAD/MCP penalized linear
regression computed by (fully convergent) LLA with the tuned LASSO
estimator as its initial solution, following \citet{Fanetal2012}.
\end{longlist}

Notice that we no longer involve LLA$_r$ and LLA$_0$ in this test,
because a similar numerical experiment presented by \citet{Fanetal2012}
has shown that LLA$_1$ is more preferable than most other LLA variants
in statistical performance.

%\end{algorithm}

%\textcolor{blue}{To Be DELETED: In the simulation, we consider two
%cases. One has $n=100$ observations and $d=500$ dimensions; the other
%one has $n=100$ observations and $d=1000$ dimensions. We also generate
%independent validation sets of a random samples of size 100 to tune
%the proposed MIPGO technique with the same validation error as in
%\citet{Fanetal2012}. }
%({\mathbf what is the number of simulations? The number
%of simulations is provided soon after the introduction of the
%alternative
%schemes})
%i.e., $\sum_{i\in validation}(y_i-x_i^\top\hat\beta)^2$.

Numerical results are presented in Table~\ref{table 3}. According to
the table, the proposed MIPGO approach estimates the (in)significant
coefficients correctly in both SCAD and MCP penalties, and provides an
improvement on the average
AD over all the other alternative schemes. %Yet the standard error of
%the absolute deviation has an undesirable increase.
%To further illustrate the comparison of the different approaches, we
%plot the 95\% confidence interval of the mean absolute deviation in
%Figure \ref{CI95}, from which we can tell that, despite the higher
%standard error, the mean absolute deviation of MIPGO significantly
%improves over the counterpart approaches for both LR-SCAD and LR-MCP.

%t5 #&#
%
\begin{table}%[h!]
\tabcolsep=0pt
\tablewidth=250pt
\caption{Comparison of statistical performance. ``Time'' stands for the computational time in seconds. The numbers in
parenthesis are the standard errors}\label{table 3}
\begin{tabular*}{\tablewidth}{@{\extracolsep{\fill}}@{}lcccc@{}}
\hline
 & \multicolumn{4}{c@{}}{$\bolds{n=100}$\textbf{,} $\bolds{d=1000}$}\\[-6pt]
 & \multicolumn{4}{c@{}}{\hrulefill}\\
\textbf{Method} & \textbf{AD} & \textbf{FP} & \textbf{FN} & \textbf{Time}\\
\hline
 {LASSO}& 2.558 & 5.700 & 0 & 2.332 \\
&(0.047) & (0.255) & (0) & (0.108) \\[3pt]
 {GM$_1$-SCAD}& 0.526 & 0.600 & 0 & 4.167 \\
&(0.017) & (0.084) & (0) & (0.254) \\[3pt]
 {GM$_1$-MCP}& 0.543 & 0.540 & 0 & 4.42 \\
&(0.018) & (0.073) & (0) & (0.874) \\[3pt]
 {GM$_2$-SCAD}& 3.816 & 18.360 & 0 & 3.968 \\
&(0.104) & (0.655) & (0) & (0.049) \\[3pt]
 {GM$_2$-MCP}& 0.548 & 0.610 & 0 &3.916 \\
&(0.019) & (0.083) & (0) & (0.143) \\[3pt]
 {SparseNet}& 1.012 & 5.850 & 0 & 2.154 \\
&(0.086) & (1.187) & (0) & (0.017)\\[3pt]
 {Ncvreg-SCAD}& 1.068 & 9.220 & 0 & 0.733 \\
&(0.061) & (0.979) & (0) &(0.007)\\[3pt]
 {Ncvreg-MCP}& 0.830 & 3.200 & 0 & 0.877 \\
&(0.045) & (0.375) & (0) &(0.009)\\[3pt]
 {SCAD-LLA$_1$}& 0.526 & 0.600 & 0 & 31.801\\
&(0.017) & (0.084) & (0) &(1.533) \\[3pt]
 {MCP-LLA$_1$}& 0.543 & 0.540 & 0 &28.695\\
&(0.018) & (0.073) & (0) & (1.473) \\[3pt]
 {MIPGO-SCAD}& 0.509 & 0 & 0 &472.673 \\
&(0.017) & (0) & (0) &(97.982) \\[3pt]
 {MIPGO-MCP} & 0.509 & 0 & 0 &361.460\\
& (0.017) & (0) & (0) &(70.683)
\\[3pt]
{Oracle} & 0.509 &  &  &\\
& (0.017) &  &  &\\
\hline
\end{tabular*}
\end{table}

%To further illustrate the comparison of the different approaches, we
%plot the 95\% confidence interval of the mean absolute deviation in
%Figure \ref{CI95}, from which we can tell that, despite the higher
%standard error, the mean absolute deviation of MIPGO significantly
%improves over the counterpart approaches for both LR-SCAD and LR-MCP.

\begin{figure}

\includegraphics{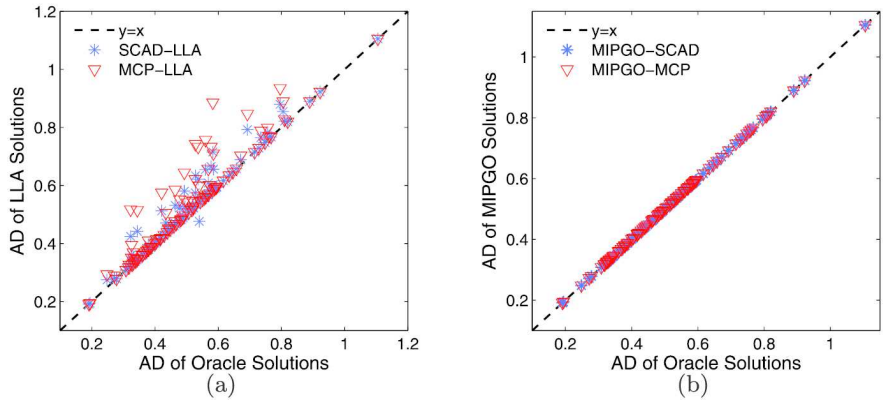}

\caption{Comparison between generated solutions and the oracle
solutions in AD when \textup{(a)} solutions are generated by LLA$_1$, and \textup{(b)}
solutions are generated by MIPGO. The horizontal axis is the AD value
of the oracle solution for each simulation, while the vertical axis is
the AD of generated solutions for the same simulation. The closer a
point is to the linear function ``$y=x$,'' the smaller is the
difference between the AD of a generated solution and the AD of the
corresponding oracle solution.}\label{compare with oracle plot}
\end{figure}

To further measure the performance of different schemes, we use the
oracle estimator as a benchmark. The oracle estimator is computed as
following: denote by $x_{t,i}$ as the $i$th dimension of the $t$th
sample $x_t$, and by $\mathcal S$ the true support set, that
is, $\mathcal S:= \{i: \beta_i^{\mathrm{true}}\neq0\}$. We conduct a linear
regression using $\hat X:= (x_{t,i}: t=1,\dots,n, i\in\mathcal S)$
and $y:= (y_t)$. As has been shown in Table~\ref{table 3}, MIPGO
yields a very close average AD and standard error to the oracle
estimator. This observation is further confirmed in Figure~\ref{compare
with oracle plot}. Specifically, Figure~\ref{compare with oracle plot}(a) and~(b) illustrate relative the
performance of LLA$_1$ and of MIPGO, respectively, in contrast to the
oracle estimators. We see that MIPGO well approximates the oracle
solution. Comparing MIPGO and LLA$_1$ from the figures, we can tell a
noticeably improved recovery quality by MIPGO in contrast to LLA$_1$.

Nonetheless, we would like to remark that, although
MIPGO yields a better solution quality over all the other local
algorithms in every cases of the experiment as presented, the local
algorithms are all
noticeably faster than MIPGO. Therefore, we think
that MIPGO is less advantageous in terms of computational time.

%This is also represented in Figures \ref{compare with oracle
%plot}(a)--(b), from which we observe that the MIPGO solutions with both
%SCAD and MCP penalties recover the oracle solutions very well in most
%of the simulations. In fact, only in two simulations with $d=1000$,
%the MIPGO-SCAD yields a relatively higher AD value, as shown in Figure
%\ref{compare with oracle plot}(b).
%
%f4 #&#
%

%s6.1 #&#
\subsection{A real data example}
In this section, we conduct our last numerical test comparing MIPGO,
LLA and the gradient methods on a real data set collected in a
marketing study
\citepp{Wang,Lanetal}, which has a total of $n=463$ daily records. For
each record,
the response variable is the number of customers and the originally
6397 predictors are sales volumes of
products. To facilitate computation, we employ the feature screening
scheme in
\citet{LiScreening} to reduce the dimension to 1500. The numerical
results are summarized in
Table~\ref{table 71}. In this table, GM$_1$ and GM$_2$ refer to the
local solution methods proposed
by \citet{LohandWainwright} and by \citet{Wangetal2013}, respectively.
LLA$_0$ denote
the LLA initialized as zero. LLA$_1$ denote the LLA initialized as the
solution generated by LASSO.
To tune the LASSO, we implement LLA$_1$ choosing the coefficients $
\omega$ in the
LASSO problem \eqref{LASSO problem} from the set $\{0.1\times nK\lambda
: K=\{0,1,\dots,20\}\}$ %({\mathbf why do this?})
and we select the $\omega$ that enables LLA$_1$ to yield the best
objective value. Here, the value of $\lambda$ is the same as the tuning
parameter of SCAD or MCP. As reported in Table~\ref{table 71}$,\lambda
=0.02$ for SCAD, and $\lambda=0.03$ for MCP, respectively. Observations
from Table~\ref{table 71} can be summarized as following: (i) for the
case with the SCAD penalty,
the proposed MIPGO yields a significantly better solution than all
other alternative schemes in terms of
both Akaike's information criterion (AIC), Bayesian information
criterion (BIC)
and the objective value. Furthermore, MIPGO also outputs a model with
the smallest number of parameters.
(ii) for the MCP case, both MIPGO and LLA$_1$ outperforms other
schemes. Yet these two approaches have
similar values for AIC and BIC. Nonetheless, MIPGO provides a better
model as the number of nonzero parameters is smaller than the solution
generated by LLA$_1$.
%t6 #&#
%
\begin{table}
\tabcolsep=0pt
\caption{Results of the Real Data Example.
``NZ''$,\#_{0.05}$ and $\#_{0.10}$ stand
for the numbers of parameters that are nonzero, that has a p-value
greater or equal to 0.05,
and that has a p-value greater or equal to 0.1. ``$R^2$'' denotes the
R-squared value. ``AIC,'' ``BIC'' and ``Obj.'' stand for Akaike's information
criterion, Bayesian information criterion and objective function value}\label{table 71}
\begin{tabular*}{\tablewidth}{@{\extracolsep{\fill}}@{}ld{3.0}d{3.0}d{3.0}cd{4.3}d{4.3}d{3.3}@{}}
\hline
\textbf{Method} &\multicolumn{1}{c}{\textbf{NZ}} & \multicolumn{1}{c}{$\bolds{\#_{0.05}}$} &  \multicolumn{1}{c}{$\bolds{\#_{0.10}}$} & \multicolumn{1}{c}{$\bolds{R^2}$} &
\multicolumn{1}{c}{\textbf{AIC}} & \multicolumn{1}{c}{\textbf{BIC}}& \textbf{Obj.}\\
\hline
\multicolumn{8}{c@{}}{SCAD: $\lambda=0.02$; $a=3.7$}\\
 {GM$_1$}& 1500 & \multicolumn{1}{c}{--} & \multicolumn{1}{c}{--} &0.997 &357.101 &6563.691 & 212.279\\
 {GM$_2$}& 401 & 401 &401 &0.698 &246.886 &1906.114 & 103.626 \\
 {LLA$_{0}$} & 185 & 119 & 115 &0.864&-554.093&211.387& 76.031\\
 {LLA$_{1}$} & 181 & 83 & 80 &0.912 &-763.718&14.789&71.673\\
 {MIPGO} & 129 & 35 & 34 &0.898 &-796.581&-262.814&68.474 \\[3pt]
\multicolumn{8}{c@{}}{MCP: $\lambda=0.03$; $a=2$}\\
 {GM$_1$}& 818 & \multicolumn{1}{c}{--} & \multicolumn{1}{c}{--} & 0.332 &1448.966 & 5129.436 &169.091\\
 {GM$_2$}& 134 & 110 & 104 & 0.735 & -296.624 & -346.474 &93.870\\
 {LLA$_{0}$} & 96 & 5 & 6 & 0.856 &704.654 &-307.432 &72.645 \\
 {LLA$_{1}$} & 113 & 2 & 2 & 0.902 &-849.842 &-382.279 &69.292 \\
 {MIPGO} & 109 & 3 & 3 & 0.899 &-841.280 &-390.267 &68.591 \\
 \hline
\end{tabular*}
\end{table}

%
%s7 #&#
\section{Conclusion}\label{conclusion}
The lack of solution schemes that ascertain solution quality to
nonconvex learning
with folded concave penalty has been an open problem in sparse
recovery. In this paper, we seek to address
this issue in a direct manner by proposing a global optimization
technique for a class of nonconvex learning
problems without imposing very restrictive conditions.
%(general in the sense that no artificial assumptions are additional
%stipulated to the original problem.)

In this paper, we provide a reformulation of the nonconvex learning
problem into a general quadratic program. This reformulation then
enables us to have the following findings:
\begin{longlist}[(a)]
\item[(a)]To formally state the complexity of finding the global
optimal solution to the nonconvex learning with the SCAD and the MCP penalties.
\item[(b)]To derive a MIP-based global optimization approach, MIPGO, to
solve the SCAD and MCP penalized nonconvex learning problems with
theoretical guarantee. Numerical results indicate that the proposed MIPGO
outperforms the gradient method by \citet{LohandWainwright} and
\citet
{Wangetal2013} and LLA approach
with different initialization schemes in solution quality and
statistical performance.
\end{longlist}

To the best of our knowledge, the complexity bound of solving the
nonconvex learning with the MCP
and SCAD penalties globally has not been reported in literature and
MIPGO is the first
optimization scheme with provable guarantee on global optimality for
solving a folded concave penalized learning problem.

We would like to alert the readers that the proposed MIPGO scheme,
though being
effective in globally solving the nonconvex learning with the MCP and SCAD
penalty problem, yields a comparatively larger computational overhead
than the
local solution method in larger scale problems. (See comparison of
computing times
in Table~\ref{table 3}.) In the practice of highly
time-sensitive statistical learning with hugh problem sizes, LLA and
other local solution schemes can work more
efficiently. However, there are important application scenarios where a
further refinement on the solution quality or even the exact global
optimum is
required. MIPGO is particularly effective in those applications, as it is
the only method that is capable of providing the refinement with theoretical
guarantee.

Finally, we would like to remark that the quadratic programming reformulation
of penalized least squares with the MCP and SCAD penalty can be further
exploited to develop convex approximation, complexity analyses and solution
schemes for finding a local solution. Those will be the future
extensions of
the presented work herein.

\section{Proofs of Theorems \texorpdfstring{\protect\ref{prop 1}}{3.2} and \texorpdfstring{\protect\ref
{LRmcp}}{3.4}}\label{some proofs} In this section, we give proofs of
Theorems \ref{prop 1} and \ref{LRmcp}.

\begin{pf*}{Proof of Theorem \ref{prop 1}}
Recall that $\mathbf1$
denotes an all-ones vector of a proper dimension. The program has a
Lagrangian $\mathbb F_{\mathrm{SCAD}}$ given as
\begin{eqnarray*}
&& \mathbb F_{\mathrm{SCAD}}(\beta,g,h,\gamma_1,\gamma_2,
\gamma_3,\gamma_4,\rho)
\\
&&\qquad := \tfrac{1}{2}\bigl[\beta^\top Q\beta+n(a-1)g^\top
g+2ng^\top h\bigr]+q^\top \beta-na\lambda{
\mathbf1}^\top g+\gamma_1^\top(\beta-h)
\\
&&\quad\qquad{} -\gamma_2^\top(\beta+h)+\gamma_3^\top(-g)+
\gamma_4^\top(g-\lambda {\mathbf1})+\rho^\top
\bigl(\mathcal A^\top\beta-\mathbf b\bigr),
\end{eqnarray*}
where $\gamma_{1}:= (\gamma_{1,i} )\in\mathbb R^d_+$, $\gamma_{2}
:= (\gamma_{2,i} )\in\mathbb R^d_+$, $\gamma_{3}:= (\gamma_{3,i} )\in
\mathbb R^d_+$, $\gamma_{4}:= (\gamma_{4,i} )\in\mathbb R^d_+$, and
$\rho\in\mathbb R_+^m$ are Lagrangian multipliers. The KKT condition yields
%e8.1 #&#
%e8.2 #&#
%
\begin{eqnarray}
&&\cases{ \displaystyle\frac{\partial\mathbb F_{\mathrm{SCAD}}}{\partial\beta}:= Q\beta+q+\gamma _1-
\gamma_2+\mathcal A\rho=0,%\label{KKT SCAD 1}
\vspace*{2pt}\cr
\displaystyle \frac{\partial\mathbb F_{\mathrm{SCAD}}}{\partial h}:= ng-
\gamma_1-\gamma_2=0,%
%\label{KKT SCAD 1a}
\vspace*{2pt}\cr
\displaystyle\frac{\partial\mathbb F_{\mathrm{SCAD}}}{\partial g}:= n(a-1)g+nh-na\lambda {\mathbf1}-\gamma_3+
\gamma_4=0,}\label{KKT SCAD 1b}
\\
\label{KKT SCAD 2} &&\left\{
\begin{array}{l}
\left.
\begin{array}{l}
\gamma_{1,i}\geq0;\qquad\gamma_{1,i}\cdot(\beta_i-h_i)=0,\\
\gamma_{2,i}\geq0;\qquad \gamma_{2,i}\cdot(-\beta_i-h_i)=0,\\
\gamma_{3,i}\geq0;\qquad\gamma_{3,i}\cdot g_i=0,\\
\gamma_{4,i}\geq0;\qquad\gamma_{4,i}\cdot(g_i-\lambda)=0,
\end{array} \right\} \qquad\forall i=1,\dots,d.
\\
\rho\geq0;\qquad \rho^\top\bigl(\mathcal A^\top\beta-\mathbf b \bigr)=0.
\end{array}\right.
\end{eqnarray}
Since $\Lambda$ is nonempty and $\mathcal A$ is full rank, it is easy
to check that the linear independence constraint qualification is
satisfied. Therefore, the global solution satisfies the KKT condition.
This leads us to an equivalent representation of \eqref{quadratic obj} %
%\eqref{quadratic const2}
in the form:
%e8.3 #&#
%e8.4 #&#
%
\begin{eqnarray}
\qquad&& \min \tfrac{1}{2}\bigl[\beta^\top Q\beta+n(a-1)g^\top
g+2ng^\top h\bigr]+q^\top \beta-na\lambda{
\mathbf1}^\top g,\qquad \mbox{s.t.}\label{add const 1}
\\
&&\qquad \cases{ \beta\in\Lambda;\qquad h\geq\beta;\qquad h\geq-\beta;\qquad0\leq g\leq
\lambda,
\vspace*{2pt}\cr
\mbox{Constraints \eqref{KKT SCAD 1b}--\eqref{KKT SCAD 2}},
\vspace*{2pt}\cr
\beta, g, h, \gamma_1, \gamma_2, \gamma_3,
\gamma_4\in\mathbb R^d_+;\qquad  \rho\in\mathbb
R^m_+.}\label{add const 1a}
\end{eqnarray}
Then it suffices to show that \eqref{add const 1}--\eqref{add const 1a}
is equivalent to \eqref{to show comple111}--\eqref{comple total 1}.
%-\eqref{comple4}.

Notice that the objective function \eqref{add const 1} is immediately
\[
\tfrac{1}{2} \beta^\top(Q\beta+q)+\tfrac{1}{2}q^\top
\beta+\tfrac
{1}{2}g^\top \bigl(n(a-1)g-na\lambda{\mathbf1}+2nh
\bigr)-\tfrac{1}{2}na\lambda{\mathbf 1}^\top g =: I_1.
%\label{obj interim 1}
\]
Due to equalities \eqref{KKT SCAD 1b}, %\eqref{KKT SCAD 1}, \eqref{KKT
%SCAD 1a},,
%
\begin{eqnarray}
I_1&=&%\frac{1}{2}\beta^\top(\gamma_2-\gamma_1-\mathcal A \rho)-
%\frac{1}{2}q^\top\beta+\frac{1}{2}g^\top(\gamma_3-\gamma_4+nh)-
%\frac{1}{2}na\lambda e^\top g
\tfrac{1}{2}
\beta^\top(\gamma_2-\gamma_1)-\tfrac{1}{2}
\beta^\top\mathcal A\rho+\tfrac{1}{2}q^\top\beta
\nonumber\\[-8pt]\\[-8pt]\nonumber
&&{}+\tfrac{1}{2}h^\top(\gamma_1+\gamma_2)+
\tfrac{1}{2}g^\top(\gamma _3-\gamma_4)-
\tfrac{1}{2}na\lambda{\mathbf1}^\top g.
\end{eqnarray}
Invoking the complementarity conditions in \eqref{KKT SCAD 2}, we may have
%e8.5 #&#
%
\begin{eqnarray}
I_1=\tfrac{1}{2}q^\top\beta-\tfrac{1}{2}\mathbf
b^\top \rho-\tfrac{1}{2}na\lambda{\mathbf1}^\top g-
\tfrac{1}{2}\lambda\gamma _4^\top {
\mathbf1}.\label{final obj SCAD}
\end{eqnarray}
Therefore, Program \eqref{add const 1}--\eqref{add const 1a} is
equivalent to
%e8.6 #&#
%
\begin{equation}
\min I_1= \tfrac{1}{2}q^\top\beta-\tfrac{1}{2}
\mathbf b^\top \rho-\tfrac{1}{2}na\lambda{\mathbf1}^\top
g-\frac{1}{2}\lambda\gamma _4^\top {\mathbf1}\qquad\mbox{s.t. \eqref{add const 1a}},
\end{equation}
which is immediately the desired result.
\end{pf*}

\begin{pf*}{Proof of Theorem \ref{LRmcp}}
The proof follows a
closely similar argument as that for Theorem \ref{LRmcp}.
The Lagrangian $\mathbb F_{\mathrm{MCP}}$ of Program \eqref{reformulated LR-MCP
obj} yields
\begin{eqnarray}
&& \mathbb F_{\mathrm{MCP}}(\beta,g,h,\eta_1,\eta_2,
\eta_3,\eta_4,\rho)\nonumber
\\
&&\qquad :=\frac{1}{2}\beta^\top Q\beta+q^\top\beta+n
\frac{1}{2a}g^\top g-n \biggl(\frac{1}{a}g- \lambda\mathbf1
\biggr)^\top h +\eta_1^\top(\beta-h)
\\
&&\quad\qquad{} +\eta_2^\top(-\beta-h)-\eta_3^\top
g+\eta_4^\top (g-a\lambda\mathbf1)+\rho^\top
\bigl(\mathcal A^\top\beta-\mathbf b\bigr),\nonumber
\end{eqnarray}
where $\mathbf1$ denotes an all-ones vector of a proper dimension, and
where $\eta_1$, $\eta_2$, $\eta_3$, $\eta_4\in\mathbb R^d_+$ and $\rho
\in\mathbb R^m_+$ are Lagrangian multipliers. The KKT condition yields
%e8.7 #&#
%e8.8 #&#
%
\begin{eqnarray}
&& \cases{\displaystyle\frac{\partial\mathbb F_{\mathrm{MCP}}}{\partial\beta}:= q+Q\beta+\eta _1 -\eta_2+
\mathcal A\rho=0,
\vspace*{2pt}\cr
\displaystyle\frac{\partial\mathbb F_{\mathrm{MCP}}}{\partial g}:= \frac{n}{a}g-\frac{n}{a}
h-\eta_3+\eta_4=0,
\vspace*{2pt}\cr
\displaystyle\frac{\partial\mathbb F_{\mathrm{MCP}}}{\partial h}:= -n \biggl(
\frac{1}{a}g- \lambda\mathbf1 \biggr)-\eta_1-
\eta_2=0,}\label{equality MCP 2}
\\
&&\cases{ \eta_1^\top(\beta-h)=0;\qquad \eta_2^\top(-
\beta-h)=0,
\vspace*{2pt}\cr
\eta_3^\top g=0;\qquad\eta_4^\top(g-a
\lambda\mathbf1)=0;\qquad\rho ^\top\bigl(\mathcal A^\top\beta-
\mathbf b\bigr)=0,
\vspace*{2pt}\cr
\eta_1\geq0;\qquad\eta_2\geq0;
\qquad\eta_3\geq0;\qquad\eta_4\geq 0;\qquad\rho\geq0.}
\label{KKT MCP 2}
\end{eqnarray}
Since $\Lambda$ is nonempty and $\mathcal A$ is full rank, it is
easily verifiable that the linear independence constraint qualification
is satisfied. This means the KKT system holds at the global solution.
Therefore, imposing additional constraints \eqref{equality MCP
2}--\eqref{KKT MCP 2} in program \eqref{reformulated LR-MCP obj}
%-\eqref{reformulated LR-MCP const}
will not result in inequivalence. Notice that the objective function in~\eqref{reformulated LR-MCP obj} equals
\begin{eqnarray}\label{obj MCP interim}
%\mbox{Eq. \eqref{reformulated LR-MCP obj}}
%\\=& q^\top\beta+\frac{1}{2}\beta^\top Q\beta+n\frac{1}{2a}g^\top g-
%\frac{n}{a}g^\top h+\lambda n\mathbf1^\top h
I_2&:=& \frac{1}{2}q^\top
\beta+\biggl(\frac{1}{2}q^\top+\frac{1}{2}
\beta^\top Q\biggr)\beta+\frac{n}{2a}g^\top(g-h)
\nonumber\\[-8pt]\\[-8pt]\nonumber
&&{} -\frac{1}{2}n \biggl(\frac{1}{a}g^\top-\lambda
\mathbf1^\top \biggr) h+\frac{1}{2}\lambda n\mathbf1^\top
h.
\end{eqnarray}
Per \eqref{equality MCP 2}, we obtain
\begin{eqnarray}\label{MCP inter 2}
I_2&=& \tfrac{1}{2}q^\top\beta-\tfrac{1}{2}(
\eta_1 -\eta_2+\mathcal A \rho )^\top\beta
\nonumber\\[-8pt]\\[-8pt]\nonumber
&&{} +\tfrac{1}{2}g^\top(\eta_3-\eta_4)+
\tfrac{1}{2}(\eta_1+\eta_2)^\top h+
\tfrac{1}{2}\lambda n\mathbf1^\top h.
\end{eqnarray}
Further noticing \eqref{KKT MCP 2}, we obtain
\begin{eqnarray*}
I_2&=& \tfrac{1}{2}q^\top\beta-\tfrac{1}{2}
\mathbf b^\top\rho+\tfrac
{1}{2}g^\top(
\eta_3-\eta_4)+\tfrac{1}{2}\lambda n
\mathbf1^\top h
\nonumber
\\
&=& \tfrac{1}{2}q^\top\beta-\tfrac{1}{2}\mathbf
b^\top\rho-\tfrac
{1}{2}g^\top\eta_4+
\tfrac{1}{2}\lambda n\mathbf1^\top h
\nonumber
\\
&=& \tfrac{1}{2}q^\top\beta-\tfrac{1}{2}\mathbf
b^\top\rho-\tfrac
{1}{2}a\lambda\mathbf1^\top
\eta_4+\tfrac{1}{2}\lambda n\mathbf1^\top h
\nonumber
\end{eqnarray*}
which immediately leads to the desired result.
\end{pf*}

%\begin{appendix}
%\section{}
%\end{appendix}

%\section*{Acknowledgments}

\begin{supplement}[id=suppA]
\stitle{Supplement to ``Global solutions to folded concave penalized nonconvex learning''}
\slink[doi]{10.1214/15-AOS1380SUPP} %[doi,text={...}] - jei reikia
%suskaldyti doi
\sdatatype{.pdf}
\sfilename{aos1380\_supp.pdf}
\sdescription{This supplemental material includes the proofs of
Proposition 2.1, 2.3 and Lemma 4.1, and some additional numerical results.}
\end{supplement}

% imsref loaded by linak, 2015-10-20 15:57:54
%
% imsref loaded by linak, 2016-01-12 10:38:13

\printaddresses

\begin{thebibliography}{30}

%b1 ###
\bibitem[\protect\citeauthoryear{Bertsimas, Chang and
Rudin}{2011}]{Bertsimasetal2011}
%
\begin{bmisc}[auto:parserefs-M02]
\bauthor{\bsnm{Bertsimas},~\bfnm{D.}\binits{D.}},
\bauthor{\bsnm{Chang},~\bfnm{A.}\binits{A.}} \AND
\bauthor{\bsnm{Rudin},~\bfnm{C.}\binits{C.}}
(\byear{2011}).
\bhowpublished{Integer optimization methods for supervised ranking.
Available at \surl{http://hdl.handle.net/1721.1/67362}}.
\end{bmisc}
%

\bptok{imsref}%
% NOT OUTPUTTED:
% url = http://hdl.handle.net/1721.1/67362
\endbibitem

%b2 ###
\bibitem[\protect\citeauthoryear{Breheny and Huang}{2011}]{BrehenyandHuang2011}
%
\begin{barticle}[mr]
\bauthor{\bsnm{Breheny},~\bfnm{Patrick}\binits{P.}} \AND
\bauthor{\bsnm{Huang},~\bfnm{Jian}\binits{J.}}
(\byear{2011}).
\btitle{Coordinate descent algorithms for nonconvex penalized
regression, with applications to biological feature selection}.
\bjournal{Ann. Appl. Stat.}
\bvolume{5}
\bpages{232--253}.
\bid{doi={10.1214/10-AOAS388}, issn={1932-6157}, mr={2810396}}
\end{barticle}
%

\bptok{imsref}%
% NOT OUTPUTTED:
% number = 1
% doi = http://dx.doi.org/10.1214/10-AOAS388
% fjournal = The Annals of Applied Statistics
\endbibitem

%b3 ###
\bibitem[\protect\citeauthoryear{Fan and Li}{2001}]{FanandLi2001}
%
\begin{barticle}[mr]
\bauthor{\bsnm{Fan},~\bfnm{Jianqing}\binits{J.}} \AND
\bauthor{\bsnm{Li},~\bfnm{Runze}\binits{R.}}
(\byear{2001}).
\btitle{Variable selection via nonconcave penalized likelihood and its
oracle properties}.
\bjournal{J. Amer. Statist. Assoc.}
\bvolume{96}
\bpages{1348--1360}.
\bid{doi={10.1198/016214501753382273}, issn={0162-1459}, mr={1946581}}
\end{barticle}
%

\bptok{imsref}%
% NOT OUTPUTTED:
% number = 456
% doi = http://dx.doi.org/10.1198/016214501753382273
% coden = JSTNAL
% fjournal = Journal of the American Statistical Association
\endbibitem

%b4 ###
\bibitem[\protect\citeauthoryear{Fan and Lv}{2011}]{FanandLv2011}
%
\begin{barticle}[mr]
\bauthor{\bsnm{Fan},~\bfnm{Jianqing}\binits{J.}} \AND
\bauthor{\bsnm{Lv},~\bfnm{Jinchi}\binits{J.}}
(\byear{2011}).
\btitle{Nonconcave penalized likelihood with NP-dimensionality}.
\bjournal{IEEE Trans. Inform. Theory}
\bvolume{57}
\bpages{5467--5484}.
\bid{doi={10.1109/TIT.2011.2158486}, issn={0018-9448}, mr={2849368}}
\end{barticle}
%

\bptok{imsref}%
% NOT OUTPUTTED:
% number = 8
% doi = http://dx.doi.org/10.1109/TIT.2011.2158486
% coden = IETTAW
% fjournal = Institute of Electrical and Electronics Engineers.
%Transactions on Information Theory
\endbibitem

%b5 ###
\bibitem[\protect\citeauthoryear{Fan, Xue and Zou}{2014}]{Fanetal2012}
%
\begin{barticle}[mr]
\bauthor{\bsnm{Fan},~\bfnm{Jianqing}\binits{J.}},
\bauthor{\bsnm{Xue},~\bfnm{Lingzhou}\binits{L.}} \AND
\bauthor{\bsnm{Zou},~\bfnm{Hui}\binits{H.}}
(\byear{2014}).
\btitle{Strong oracle optimality of folded concave penalized estimation}.
\bjournal{Ann. Statist.}
\bvolume{42}
\bpages{819--849}.
\bid{doi={10.1214/13-AOS1198}, issn={0090-5364}, mr={3210988}}
\end{barticle}
%

\bptok{imsref}%
% NOT OUTPUTTED:
% number = 3
% doi = http://dx.doi.org/10.1214/13-AOS1198
% fjournal = The Annals of Statistics
\endbibitem

%b6 ###
\bibitem[\protect\citeauthoryear{Grant and Boyd}{2008}]{CVX2}
%
\begin{bincollection}[mr]
\bauthor{\bsnm{Grant},~\bfnm{Michael~C.}\binits{M.~C.}} \AND
\bauthor{\bsnm{Boyd},~\bfnm{Stephen~P.}\binits{S.~P.}}
(\byear{2008}).
\btitle{Graph implementations for nonsmooth convex programs}.
In \bbooktitle{Recent Advances in Learning and Control}.
\bseries{Lecture Notes in Control and Inform. Sci.}
\bvolume{371}
\bpages{95--110}.
\bpublisher{Springer},
\blocation{London}.
\bid{doi={10.1007/978-1-84800-155-8_7}, mr={2409077}}
\end{bincollection}\vadjust{\goodbreak}
%

\bptok{imsref}%
% NOT OUTPUTTED:
% doi = http://dx.doi.org/10.1007/978-1-84800-155-8_7
\endbibitem

%b7 ###
\bibitem[\protect\citeauthoryear{Grant and Boyd}{2013}]{CVX1}
%
\begin{bmisc}[auto:parserefs-M02]
\bauthor{\bsnm{Grant},~\bfnm{M.}\binits{M.}} \AND
\bauthor{\bsnm{Boyd},~\bfnm{S.}\binits{S.}}
(\byear{2013}).
\bhowpublished{CVX: Matlab software for disciplined convex programming,
version 2.0 beta.
Available at \surl{http://cvxr.com/cvx}}.
\end{bmisc}
%

\bptok{imsref}%
% NOT OUTPUTTED:
% url = http://cvxr.com/cvx
\endbibitem

%b8 ###
\bibitem[\protect\citeauthoryear{Huang and Zhang}{2012}]{HuangandZhang2012}
%
\begin{barticle}[mr]
\bauthor{\bsnm{Huang},~\bfnm{Jian}\binits{J.}} \AND
\bauthor{\bsnm{Zhang},~\bfnm{Cun-Hui}\binits{C.-H.}}
(\byear{2012}).
\btitle{Estimation and selection via absolute penalized convex
minimization and its multistage adaptive applications}.
\bjournal{J. Mach. Learn. Res.}
\bvolume{13}
\bpages{1839--1864}.
\bid{issn={1532-4435}, mr={2956344}}
\end{barticle}
%

\bptok{imsref}%
% NOT OUTPUTTED:
% fjournal = Journal of Machine Learning Research (JMLR)
\endbibitem

%b9 ###
\bibitem[\protect\citeauthoryear{Hunter and Li}{2005}]{HunterandLi2005}
%
\begin{barticle}[mr]
\bauthor{\bsnm{Hunter},~\bfnm{David~R.}\binits{D.~R.}} \AND
\bauthor{\bsnm{Li},~\bfnm{Runze}\binits{R.}}
(\byear{2005}).
\btitle{Variable selection using MM algorithms}.
\bjournal{Ann. Statist.}
\bvolume{33}
\bpages{1617--1642}.
\bid{doi={10.1214/009053605000000200}, issn={0090-5364}, mr={2166557}}
\end{barticle}
%

\bptok{imsref}%
% NOT OUTPUTTED:
% number = 4
% doi = http://dx.doi.org/10.1214/009053605000000200
% coden = ASTSC7
% fjournal = The Annals of Statistics
\endbibitem

%b10 ###
\bibitem[\protect\citeauthoryear{Kim, Choi and Oh}{2008}]{KimCCCP}
%
\begin{barticle}[mr]
\bauthor{\bsnm{Kim},~\bfnm{Yongdai}\binits{Y.}},
\bauthor{\bsnm{Choi},~\bfnm{Hosik}\binits{H.}} \AND
\bauthor{\bsnm{Oh},~\bfnm{Hee-Seok}\binits{H.-S.}}
(\byear{2008}).
\btitle{Smoothly clipped absolute deviation on high dimensions}.
\bjournal{J. Amer. Statist. Assoc.}
\bvolume{103}
\bpages{1665--1673}.
\bid{doi={10.1198/016214508000001066}, issn={0162-1459}, mr={2510294}}
\end{barticle}
%

\bptok{imsref}%
% NOT OUTPUTTED:
% number = 484
% doi = http://dx.doi.org/10.1198/016214508000001066
% coden = JSTNAL
% fjournal = Journal of the American Statistical Association
\endbibitem

%b11 ###
\bibitem[\protect\citeauthoryear{Lan et~al.}{2013}]{Lanetal}
%
\begin{barticle}[auto:parserefs-M02]
\bauthor{\bsnm{Lan},~\bfnm{W.}\binits{W.}},
\bauthor{\bsnm{Zhong},~\bfnm{P.-S.}\binits{P.-S.}},
\bauthor{\bsnm{Li},~\bfnm{R.}\binits{R.}},
\bauthor{\bsnm{Wang},~\bfnm{H.}\binits{H.}} \AND
\bauthor{\bsnm{Tsai},~\bfnm{C.-L.}\binits{C.-L.}}
(\byear{2013}).
\btitle{Testing a single regression coefficient in high dimensional
linear models}.
\bnote{Working paper}.
\end{barticle}
%

\bptok{imsref}%
\endbibitem

%b12 ###
\bibitem[\protect\citeauthoryear{Lawler and Wood}{1966}]{LawlerandWood}
%
\begin{barticle}[mr]
\bauthor{\bsnm{Lawler},~\bfnm{E.~L.}\binits{E.~L.}} \AND
\bauthor{\bsnm{Wood},~\bfnm{D.~E.}\binits{D.~E.}}
(\byear{1966}).
\btitle{Branch-and-bound methods: A survey}.
\bjournal{Oper. Res.}
\bvolume{14}
\bpages{699--719}.
\bid{issn={0030-364X}, mr={0202469}}
\end{barticle}
%

\bptok{imsref}%
% NOT OUTPUTTED:
% fjournal = Operations Research
\endbibitem

%b13 ###
\bibitem[\protect\citeauthoryear{Li, Zhong and Zhu}{2012}]{LiScreening}
%
\begin{barticle}[mr]
\bauthor{\bsnm{Li},~\bfnm{Runze}\binits{R.}},
\bauthor{\bsnm{Zhong},~\bfnm{Wei}\binits{W.}} \AND
\bauthor{\bsnm{Zhu},~\bfnm{Liping}\binits{L.}}
(\byear{2012}).
\btitle{Feature screening via distance correlation learning}.
\bjournal{J. Amer. Statist. Assoc.}
\bvolume{107}
\bpages{1129--1139}.
\bid{doi={10.1080/01621459.2012.695654}, issn={0162-1459}, mr={3010900}}
\end{barticle}
%

\bptok{imsref}%
% NOT OUTPUTTED:
% number = 499
% doi = http://dx.doi.org/10.1080/01621459.2012.695654
% coden = JSTNAL
% fjournal = Journal of the American Statistical Association
\endbibitem

%b14 ###
\bibitem[\protect\citeauthoryear{Liu, Yao and Li}{2016}]{LiuYaoandLi2014}
%
\begin{bmisc}[author]
\bauthor{\bsnm{Liu},~\binits{H.}},
\bauthor{\bsnm{Yao},~\binits{T.}} \AND
\bauthor{\bsnm{Li}~\binits{Runze}}
(\byear{2016}).
\bhowpublished{Supplement to ``Global solutions to folded concave
penalized nonconvex learning.''
DOI:\doiurl{10.1214/15-AOS1380SUPP}}.
\end{bmisc}
%

\bptok{imsref}%
\endbibitem

%b15 ###
\bibitem[\protect\citeauthoryear{Loh and Wainwright}{2015}]{LohandWainwright}
%
\begin{barticle}[mr]
\bauthor{\bsnm{Loh},~\bfnm{Po-Ling}\binits{P.-L.}} \AND
\bauthor{\bsnm{Wainwright},~\bfnm{Martin~J.}\binits{M.~J.}}
(\byear{2015}).
\btitle{Regularized {$M$}-estimators with nonconvexity: Statistical and
algorithmic theory for local optima}.
\bjournal{J. Mach. Learn. Res.}
\bvolume{16}
\bpages{559--616}.
\bid{issn={1532-4435}, mr={3335800}}
\bptnote{check volume, check pages}%
\end{barticle}
%

\bptok{imsref}%
% NOT OUTPUTTED:
% fjournal = Journal of Machine Learning Research (JMLR)
\endbibitem

%b16 ###
\bibitem[\protect\citeauthoryear{Mart{\'{\i}} and Reinelt}{2011}]{marti}
%
\begin{bincollection}[mr]
\bauthor{\bsnm{Mart{\'{\i}}},~\bfnm{Rafael}\binits{R.}} \AND
\bauthor{\bsnm{Reinelt},~\bfnm{Gerhard}\binits{G.}}
(\byear{2011}).
\btitle{Branch-and-bound}.
In \bbooktitle{The Linear Ordering Problem}.
\bpublisher{Springer},
\blocation{Heidelberg}.
\bid{doi={10.1007/978-3-642-16729-4}, mr={2722086}}
\end{bincollection}
%

\bptok{imsref}%
% NOT OUTPUTTED:
% doi = http://dx.doi.org/10.1007/978-3-642-16729-4
% isbn = 978-3-642-16728-7
% fpage = xii+171
\endbibitem

%b17 ###
\bibitem[\protect\citeauthoryear{Mazumder, Friedman and
Hastie}{2011}]{Mazumder2011}
%
\begin{barticle}[mr]
\bauthor{\bsnm{Mazumder},~\bfnm{R.}\binits{R.}},
\bauthor{\bsnm{Friedman},~\bfnm{J.}\binits{J.}} \AND
\bauthor{\bsnm{Hastie},~\bfnm{T.}\binits{T.}}
(\byear{2011}).
\btitle{SparseNet:  Coordinate descent with non-convex penalties}.
\bjournal{J. Amer. Statist. Assoc.}
\bvolume{106}
\bpages{1125--1138}.
\bid{mr={2894769}}
\end{barticle}
%

\bptok{imsref}%
\endbibitem

%b18 ###
\bibitem[\protect\citeauthoryear{Meinshausen and B{\"
u}hlmann}{2006}]{Meinshausen}
%
\begin{barticle}[mr]
\bauthor{\bsnm{Meinshausen},~\bfnm{Nicolai}\binits{N.}} \AND
\bauthor{\bsnm{B{\"u}hlmann},~\bfnm{Peter}\binits{P.}}
(\byear{2006}).
\btitle{High-dimensional graphs and variable selection with the lasso}.
\bjournal{Ann. Statist.}
\bvolume{34}
\bpages{1436--1462}.
\bid{doi={10.1214/009053606000000281}, issn={0090-5364}, mr={2278363}}
\end{barticle}
%

\bptok{imsref}%
% NOT OUTPUTTED:
% number = 3
% doi = http://dx.doi.org/10.1214/009053606000000281
% coden = ASTSC7
% fjournal = The Annals of Statistics
\endbibitem

%b19 ###
\bibitem[\protect\citeauthoryear{Nesterov}{2007}]{Nesterov2007}
%
\begin{bmisc}[auto:parserefs-M02]
\bauthor{\bsnm{Nesterov},~\bfnm{Y.}\binits{Y.}}
(\byear{2007}).
\bhowpublished{Gradient methods for minimizing composite objective function.
CORE Discussion Papers 2007076, Universit Catholique de Louvain, Center
for Operations Research and Dconometrics (CORE)}.
\end{bmisc}
%

\bptok{imsref}%
% NOT OUTPUTTED:
% institution = Universit Catholique de Louvain, Center for Operations
%Research and Dconometrics (CORE)
\endbibitem

%b20 ###
\bibitem[\protect\citeauthoryear{Pardalos}{1991}]{Pardalos}
%
\begin{barticle}[mr]
\bauthor{\bsnm{Pardalos},~\bfnm{Panos~M.}\binits{P.~M.}}
(\byear{1991}).
\btitle{Global optimization algorithms for linearly constrained
indefinite quadratic problems}.
\bjournal{Comput. Math. Appl.}
\bvolume{21}
\bpages{87--97}.
\bid{doi={10.1016/0898-1221(91)90163-X}, issn={0898-1221}, mr={1096136}}
\end{barticle}
%

\bptok{imsref}%
% NOT OUTPUTTED:
% number = 6-7
% doi = http://dx.doi.org/10.1016/0898-1221(91)90163-X
% coden = CMAPDK
% fjournal = Computers \& Mathematics with Applications. An
%International Journal
\endbibitem

%b21 ###
\bibitem[\protect\citeauthoryear{Vandenbussche and
Nemhauser}{2005}]{MIPliterature}
%
\begin{barticle}[mr]
\bauthor{\bsnm{Vandenbussche},~\bfnm{Dieter}\binits{D.}} \AND
\bauthor{\bsnm{Nemhauser},~\bfnm{George~L.}\binits{G.~L.}}
(\byear{2005}).
\btitle{A polyhedral study of nonconvex quadratic programs with box
constraints}.
\bjournal{Math. Program.}
\bvolume{102}
\bpages{531--557}.
\bid{doi={10.1007/s10107-004-0549-0}, issn={0025-5610}, mr={2136226}}
\end{barticle}
%

\bptok{imsref}%
% NOT OUTPUTTED:
% number = 3, Ser. A
% doi = http://dx.doi.org/10.1007/s10107-004-0549-0
% fjournal = Mathematical Programming. A Publication of the
%Mathematical Programming Society
\endbibitem

%b22 ###
\bibitem[\protect\citeauthoryear{Vavasis}{1991}]{VavasisBook}
%
\begin{bbook}[mr]
\bauthor{\bsnm{Vavasis},~\bfnm{Stephen~A.}\binits{S.~A.}}
(\byear{1991}).
\btitle{Nonlinear Optimization}.
\bseries{International Series of Monographs on Computer Science}
\bvolume{8}.
\bpublisher{The Clarendon Press, Oxford Univ. Press},
\blocation{New York}.
\bid{mr={1296253}}
\end{bbook}
%

\bptok{imsref}%
% NOT OUTPUTTED:
% isbn = 0-19-507208-1
% fpage = xiv+165
\endbibitem

%b23 ###
\bibitem[\protect\citeauthoryear{Vavasis}{1992}]{vavasis}
%
\begin{barticle}[mr]
\bauthor{\bsnm{Vavasis},~\bfnm{Stephen~A.}\binits{S.~A.}}
(\byear{1992}).
\btitle{Approximation algorithms for indefinite quadratic programming}.
\bjournal{Math. Program.}
\bvolume{57}
\bpages{279--311}.
\bid{doi={10.1007/BF01581085}, issn={0025-5610}, mr={1195028}}
\end{barticle}
%

\bptok{imsref}%
% NOT OUTPUTTED:
% number = 2, Ser. B
% doi = http://dx.doi.org/10.1007/BF01581085
% coden = MHPGA4
% fjournal = Mathematical Programming
\endbibitem

%b24 ###
\bibitem[\protect\citeauthoryear{Wang}{2009}]{Wang}
%
\begin{barticle}[mr]
\bauthor{\bsnm{Wang},~\bfnm{Hansheng}\binits{H.}}
(\byear{2009}).
\btitle{Forward regression for ultra-high dimensional variable screening}.
\bjournal{J. Amer. Statist. Assoc.}
\bvolume{104}
\bpages{1512--1524}.
\bid{doi={10.1198/jasa.2008.tm08516}, issn={0162-1459}, mr={2750576}}
\end{barticle}
%

\bptok{imsref}%
% NOT OUTPUTTED:
% number = 488
% doi = http://dx.doi.org/10.1198/jasa.2008.tm08516
% coden = JSTNAL
% fjournal = Journal of the American Statistical Association
\endbibitem

%b25 ###
\bibitem[\protect\citeauthoryear{Wang, Kim and Li}{2013}]{Wangetal2013a}
%
\begin{barticle}[mr]
\bauthor{\bsnm{Wang},~\bfnm{Lan}\binits{L.}},
\bauthor{\bsnm{Kim},~\bfnm{Yongdai}\binits{Y.}} \AND
\bauthor{\bsnm{Li},~\bfnm{Runze}\binits{R.}}
(\byear{2013}).
\btitle{Calibrating nonconvex penalized regression in ultra-high dimension}.
\bjournal{Ann. Statist.}
\bvolume{41}
\bpages{2505--2536}.
\bid{doi={10.1214/13-AOS1159}, issn={0090-5364}, mr={3127873}}
\end{barticle}
%

\bptok{imsref}%
% NOT OUTPUTTED:
% number = 5
% doi = http://dx.doi.org/10.1214/13-AOS1159
% fjournal = The Annals of Statistics
\endbibitem

%b26 ###
\bibitem[\protect\citeauthoryear{Wang, Liu and Zhang}{2014}]{Wangetal2013}
%
\begin{barticle}[mr]
\bauthor{\bsnm{Wang},~\bfnm{Zhaoran}\binits{Z.}},
\bauthor{\bsnm{Liu},~\bfnm{Han}\binits{H.}} \AND
\bauthor{\bsnm{Zhang},~\bfnm{Tong}\binits{T.}}
(\byear{2014}).
\btitle{Optimal computational and statistical rates of convergence for
sparse nonconvex learning problems}.
\bjournal{Ann. Statist.}
\bvolume{42}
\bpages{2164--2201}.
\bid{doi={10.1214/14-AOS1238}, issn={0090-5364}, mr={3269977}}
\end{barticle}
%

\bptok{imsref}%
% NOT OUTPUTTED:
% number = 6
% doi = http://dx.doi.org/10.1214/14-AOS1238
% fjournal = The Annals of Statistics
\endbibitem

%b27 ###
\bibitem[\protect\citeauthoryear{Zhang}{2010}]{Zhang2010}
%
\begin{barticle}[mr]
\bauthor{\bsnm{Zhang},~\bfnm{Cun-Hui}\binits{C.-H.}}
(\byear{2010}).
\btitle{Nearly unbiased variable selection under minimax concave penalty}.
\bjournal{Ann. Statist.}
\bvolume{38}
\bpages{894--942}.
\bid{doi={10.1214/09-AOS729}, issn={0090-5364}, mr={2604701}}
\bptnote{check volume}%
\end{barticle}
%

\bptok{imsref}%
% NOT OUTPUTTED:
% number = 2
% doi = http://dx.doi.org/10.1214/09-AOS729
% coden = ASTSC7
% fjournal = The Annals of Statistics
\endbibitem

%b28 ###
\bibitem[\protect\citeauthoryear{Zhang and Zhang}{2012}]{ZhangandZhang2012}
%
\begin{barticle}[mr]
\bauthor{\bsnm{Zhang},~\bfnm{Cun-Hui}\binits{C.-H.}} \AND
\bauthor{\bsnm{Zhang},~\bfnm{Tong}\binits{T.}}
(\byear{2012}).
\btitle{A general theory of concave regularization for high-dimensional
sparse estimation problems}.
\bjournal{Statist. Sci.}
\bvolume{27}
\bpages{576--593}.
\bid{doi={10.1214/12-STS399}, issn={0883-4237}, mr={3025135}}
\bptnote{check volume, check pages}%
\end{barticle}
%

\bptok{imsref}%
% NOT OUTPUTTED:
% number = 4
% doi = http://dx.doi.org/10.1214/12-STS399
% fjournal = Statistical Science. A Review Journal of the Institute of
%Mathematical Statistics
\endbibitem

%b29 ###
\bibitem[\protect\citeauthoryear{Zou}{2006}]{Zou2006}
%
\begin{barticle}[mr]
\bauthor{\bsnm{Zou},~\bfnm{Hui}\binits{H.}}
(\byear{2006}).
\btitle{The adaptive lasso and its oracle properties}.
\bjournal{J. Amer. Statist. Assoc.}
\bvolume{101}
\bpages{1418--1429}.
\bid{doi={10.1198/016214506000000735}, issn={0162-1459}, mr={2279469}}
\end{barticle}
%

\bptok{imsref}%
% NOT OUTPUTTED:
% number = 476
% doi = http://dx.doi.org/10.1198/016214506000000735
% coden = JSTNAL
% fjournal = Journal of the American Statistical Association
\endbibitem

%b30 ###
\bibitem[\protect\citeauthoryear{Zou and Li}{2008}]{ZouandLi2008}
%
\begin{barticle}[mr]
\bauthor{\bsnm{Zou},~\bfnm{Hui}\binits{H.}} \AND
\bauthor{\bsnm{Li},~\bfnm{Runze}\binits{R.}}
(\byear{2008}).
\btitle{One-step sparse estimates in nonconcave penalized likelihood models}.
\bjournal{Ann. Statist.}
\bvolume{36}
\bpages{1509--1533}.
\bid{doi={10.1214/009053607000000802}, issn={0090-5364}, mr={2435443}}
\end{barticle}
%

\bptok{imsref}%
% NOT OUTPUTTED:
% number = 4
% doi = http://dx.doi.org/10.1214/009053607000000802
% coden = ASTSC7
% fjournal = The Annals of Statistics
\endbibitem

\end{thebibliography}
\end{document}